\documentclass{gtart}

\input gtmonout
\volumenumber{2}
\volumeyear{1999}
\volumename{Proceedings of the Kirbyfest}
\pagenumbers{233}{258}
\papernumber{13}
\received{13 October 1999}
\published{19 November 1999}

\usepackage{amsfonts}
\usepackage{amsmath}
\usepackage{amsthm}
\usepackage{graphicx}

\DeclareMathSymbol\rightsquigarrow  {\mathrel}{AMSa}{"20}
\let\Bbb\mathbb
\def\S{section~}

\newtheorem{theorem}{Theorem}

\theoremstyle{definition}

\numberwithin{equation}{section}

\newcommand{\head}[1]{\smallbreak\noindent\bf
#1\rm\smallskip\nobreak}
\newcommand{\ub}[1]{{\bf #1}\qua}
\newcommand{\foot}{\setcounter{footnote}{1}\footnote}

\def\Z{\Bbb Z}

\def\R{\Bbb R}
\def\C{\Bbb C}

\newcommand{\e}{\text{E$_8$}}
\newcommand{\res}{\mathrm{res}}

\newcommand{\place}[3]
 {\text{\kern#1pt
    \smash{\raise#2pt\hbox{\rlap{#3}}}
  \kern-#1pt\kern-1ex}}
\newcommand{\fig}[2]{\includegraphics[scale=#2]{#1.eps}}
\newcommand{\figno}[1]
 {\medskip\centerline{\small Figure #1}}
\newcommand{\tabno}[2]
 {\medskip\centerline{\small Table #1: #2}}

\def\zero{\begin{picture}(16,14)(0,8)
  \multiput(0,0)(16,0){2}{\multiput(0,0)(0,16){2}{\circle*{3}}}
  \put(8,8){\circle*{3}}
  \put(0,0){\line(1,1){16}}
  \put(0,16){\line(1,-1){16}}
  \multiput(-6,-7)(24,0){2}{\multiput(0,0)(0,16){2}{$^1$}}
  \put(6,8){$^2$}
  \end{picture}}
\def\bee{\begin{picture}(64,16)(0,10)
  \multiput(0,0)(64,0){2}{\multiput(0,0)(0,20){2}{\circle*{3}}}
  \multiput(10,10)(44,0){2}{\circle*{3}}
  \multiput(20,10)(24,0){2}{\circle*{3}}
  \multiput(10,10)(30,0){2}{\line(1,0){14}}
  \multiput(28,10)(4,0){3}{\circle*{1}}
  \put(0,0){\line(1,1){10}}
  \put(0,20){\line(1,-1){10}}
  \put(64,0){\line(-1,1){10}}
  \put(64,20){\line(-1,-1){10}}
  \multiput(-6,-7)(72,0){2}{\multiput(0,0)(0,21){2}{$^1$}}
  \multiput(8,9)(34,0){2}{\multiput(0,0)(10,0){2}{$^2$}}
  \end{picture}}
\def\two{\begin{picture}(80,10)(-5,5)
  \multiput(0,10)(10,0){8}{\circle*{3}}
  \multiput(0,10)(10,0){7}{\line(1,0){10}}
  \put(20,0){\circle*{3}}
  \put(20,0){\line(0,1){10}}
  \put(-2,9){$^2$}
  \put(8,9){$^4$}
  \put(18,9){$^6$}
  \put(28,9){$^5$}
  \put(38,9){$^4$}
  \put(48,9){$^3$}
  \put(58,9){$^2$}
  \put(68,9){$^1$}
  \put(22,-7){$^3$}
  \end{picture}}
\def\three{\begin{picture}(60,10)(0,5)
  \multiput(0,10)(10,0){7}{\circle*{3}}
  \multiput(0,10)(10,0){6}{\line(1,0){10}}
  \put(30,0){\circle*{3}}
  \put(30,0){\line(0,1){10}}
  \multiput(-2,9)(60,0){2}{$^1$}
  \multiput(8,9)(40,0){2}{$^2$}
  \multiput(18,9)(20,0){2}{$^3$}
  \put(28,9){$^4$}
  \put(32,-7){$^2$}
  \end{picture}}
\def\four{\begin{picture}(40,20)(0,10)
  \multiput(0,20)(10,0){5}{\circle*{3}}
  \multiput(0,20)(10,0){4}{\line(1,0){10}}
  \multiput(20,0)(0,10){2}{\circle*{3}}
  \multiput(20,0)(0,10){2}{\line(0,1){10}}
  \multiput(-2,19)(40,0){2}{$^1$}
  \multiput(8,19)(20,0){2}{$^2$}
  \put(18,19){$^3$}
  \put(22,4){$^2$}
  \put(22,-7){$^1$}
  \end{picture}}

\def\zerog{\begin{picture}(36,24)(0,-2)
  \setlength{\unitlength}{.8pt}
  \put(25,8.66){\circle{10}}
  \put(10,0){\circle{10}}
  \put(15,0){\line(1,0){10}}
  \multiput(30,0)(10,0){2}{\circle{10}}
  \put(25,-8.66){\circle{10}}
  \end{picture}}
\def\beeg{\begin{picture}(60,24)(0,-2)
  \setlength{\unitlength}{.8pt}
  \multiput(10,8.66)(50,0){2}{\circle{10}}
  \multiput(15,0)(40,0){2}{\circle{10}}
  \put(25,0){\oval(10,10)[l]}
  \multiput(30,0)(5,0){3}{\circle*{1}}
  \put(45,0){\oval(10,10)[r]}
  \multiput(10,-8.66)(50,0){2}{\circle{10}}
  \end{picture}}
\def\twog{\begin{picture}(64,30)(0,-7.5)
  \setlength{\unitlength}{.8pt}
  \multiput(15,13)(0,-17.33){2}{\circle{10}}
  \multiput(20,4.33)(10,0){6}{\circle{10}}
  \put(10,-13){\circle{10}}
  \end{picture}}
\def\threeg{\begin{picture}(52,36)(0,-8)
  \setlength{\unitlength}{.8pt}
  \multiput(10,17.33)(5,-8.66){3}{\circle{10}}
  \multiput(25,-8.66)(10,0){4}{\circle{10}}
  \put(20,-17.33){\circle{10}}
  \end{picture}}
\def\fourg{\begin{picture}(40,36)(0,-8)
  \setlength{\unitlength}{.8pt}
  \multiput(35,8.66)(5,8.66){2}{\circle{10}}
  \multiput(10,0)(10,0){3}{\circle{10}}
  \multiput(35,-8.66)(5,-8.66){2}{\circle{10}}
  \end{picture}}

\begin{document}

\title{The $E_8$--manifold, singular fibers and\\handlebody decompositions}  
\asciititle{The E_8-manifold, singular fibers and 
handlebody decompositions}

\authors{Robion Kirby\\Paul Melvin}
\address{University of California, Berkeley, CA 94720, USA \\ 
Bryn Mawr College, Bryn Mawr, PA 19010, USA}
\email{kirby@math.berkeley.edu, pmelvin@brynmawr.edu} 
\begin{abstract}
The $E_8$--manifold has several natural framed link descriptions, and we
give an efficient method (via ``grapes'') for showing that they are indeed
the same 4--manifold.  This leads to explicit handle pictures for the
perturbation of singular fibers in an elliptic surface to a collection of
fishtails.  In the same vein, we show how the degeneration of a regular
fiber to a singular fiber in an elliptic surface provides rich examples of
Gromov's compactness theorem.
\end{abstract}
\asciiabstract{The E_8-manifold has several natural framed link 
descriptions, and we give an efficient method (via `grapes') for
showing that they are indeed the same 4-manifold.  This leads to
explicit handle pictures for the perturbation of singular fibers in an
elliptic surface to a collection of fishtails.  In the same vein, we
show how the degeneration of a regular fiber to a singular fiber in an
elliptic surface provides rich examples of Gromov's compactness
theorem.}

\primaryclass{57N13}

\secondaryclass{57R65,14J27}


\maketitle


\setcounter{section}{-1}
\section{Introduction}

The {\it \e--manifold\/} is the $4$--manifold obtained by plumbing together
eight copies of the cotangent disk bundle of the 2--sphere according to
the Dynkin diagram for the exceptional Lie group \e\ (Figure\ 0.1a).  As a
handlebody, this is given by the framed link shown in Figure
0.1b \cite{Kir}.\foot{Sometimes the tangent bundle is used, giving $+2$
framings, but by changing the orientation of the 4--manifold, which negates
the linking matrix of the corresponding framed link, and then the
orientation of alternate $2$--spheres to restore the off diagonal elements in
the linking matrix, we get an orientation reversing diffeomorphism between
these two descriptions;  as a complex manifold, $-2$ is natural.}  The
boundary of \e\ is the Poincar\'e homology sphere (see for example
\cite{KS}).


\begin{figure}[ht!]
\centerline{\kern10pt \raise15pt\hbox{{\fig{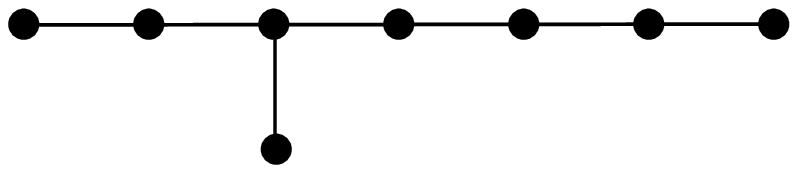}{.5}}} \kern 40pt
\fig{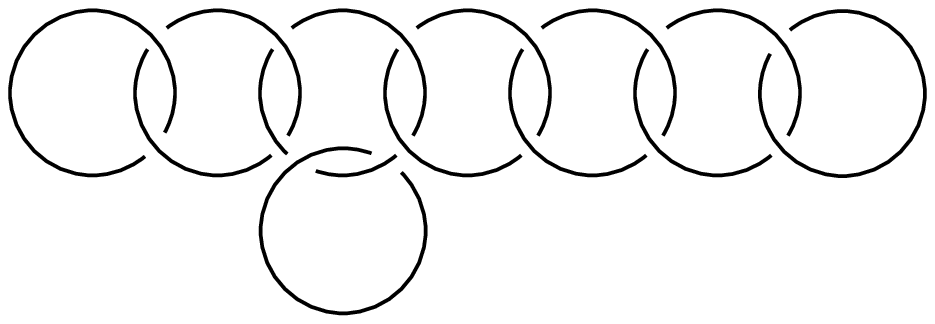}{.6}
\place{-160}{51}{$^{-2}$}
\place{-140}{51}{$^{-2}$}
\place{-120}{51}{$^{-2}$}
\place{-100}{51}{$^{-2}$}
\place{-80}{51}{$^{-2}$}
\place{-60}{51}{$^{-2}$}
\place{-40}{51}{$^{-2}$}
\place{-108}{0}{$^{-2}$}}
\medskip
\centerline{\small (a) \e--plumbing \kern 100pt (b) \e--link \kern40pt}
\figno{0.1}
\end{figure}

Alternatively, \e\ may be obtained by taking the $p$--fold cover of the
$4$--ball branched over the standard Seifert surface for the $(q,r)$--torus
knot (pushed into the interior of $B^4$) where $(p,q,r)$ is a cyclic
permutation of $(2,3,5)$.

In section 1, the calculus of framed links \cite{Kir} is used to prove that
these four 4--manifolds (\e\ and the three branched covers) are
diffeomorphic.  This result is not new.  Algebraic geometers knew this at
least as long ago as Kodaira, and it is a special case of work of Brieskorn
\cite{Br1,Br2} which we outline now.

Consider the solution $V_\varepsilon$ to $x^2 + y^3 + z^5 = \varepsilon$
in $B^6 \subset \C^3$.  This variety is a non-singular 4--manifold (for
small $\varepsilon \neq 0$) which can be described as, for example, the
$2$--fold branched cover of $B^4$ along the curve $y^3 + z^5 = \varepsilon$
(well known to be the usual Seifert surface for the $(3,5)$--torus knot). 
Similarly $V_\varepsilon$ can be viewed as a $3$ or $5$--fold branched cover.

The variety $V_0$, equal to  $V_\varepsilon$ for $\varepsilon = 0$, is a
cone on $\partial V_\varepsilon$ and has an isolated singularity at the
origin. The singularity can be resolved to obtain a non-singular complex
surface, called $V_\res$ (see \cite{HKK} for an exposition for topologists
of resolving singularities).   Brieskorn proved that $V_\res$ is
diffeomorphic to  $V_\varepsilon$ when the isolated singular point is a
{\it simple\/} singularity or a {\it rational double point}, and these are
related to the simple Lie algebras \cite{Br3}.  If these $4$--manifolds are
described using framed links, then the algebraic--geometrical proofs do not
immedately give a procedure for passing from one framed link to the other;
in particular it is not clear how complicated such a procedure might be.  So
a method is given in section 1.  The steps in Figure 1.7 from the \e--link
to the ``bunch of grapes" (Figure\ 1.3b)  are the most interesting.

Section 2 of the paper is concerned with the various singular fibers
that can occur in an elliptic surface.  These were classified by
Kodaira \cite{Kod} and a description for topologists can be found in
\cite{HKK} or the book of Gompf and Stipsicz \cite{GS} (also see \S2).  A
singular fiber, when perturbed, breaks up into a finite number of the
simplest singular fibers;  these are called {\it fishtails\/} and each
consists of an immersed $2$--sphere with one double point.  Thus a
neighborhood of a singular fiber should be diffeomorphic to a neighborhood
of several fishtails, and this is known to be diffeomorphic to a thickened
regular fiber, $T^2 \times B^2$, with several $2$--handles attached to
vanishing cycles.  Constructing these diffeomorphisms is the subject of
Section 2.

This can be looked at from a different perspective.  Gromov's compactness
theorem \cite{Gro} for (pseudo)holomorphic curves in (almost) complex
surfaces says that a sequence of curves can only degenerate in the limit by
pinching loops in the domain so as to bubble off $2$--spheres, and then
mapping the result by a holomorphic map (often just a branched covering)
onto its image.  When showing that fishtails equal a singular fiber, one
gets an idea of how a torus fiber degenerates onto the singular fiber (the
limiting curve in Gromov's sense).  Section 3 contains a discussion of this
degeneration for each of Kodaira's singular fibers.


\section{Handlebody descriptions of the \e--manifold}

The \e--manifold can be constructed by adding 2--handles to the framed
link in $S^3 = \partial B^4$ drawn in Figure 0.1b, and again in Figures
1.1a and 1.1b using successively abbreviated notation which we now explain.


\begin{figure}[ht!]
\centerline{\fig{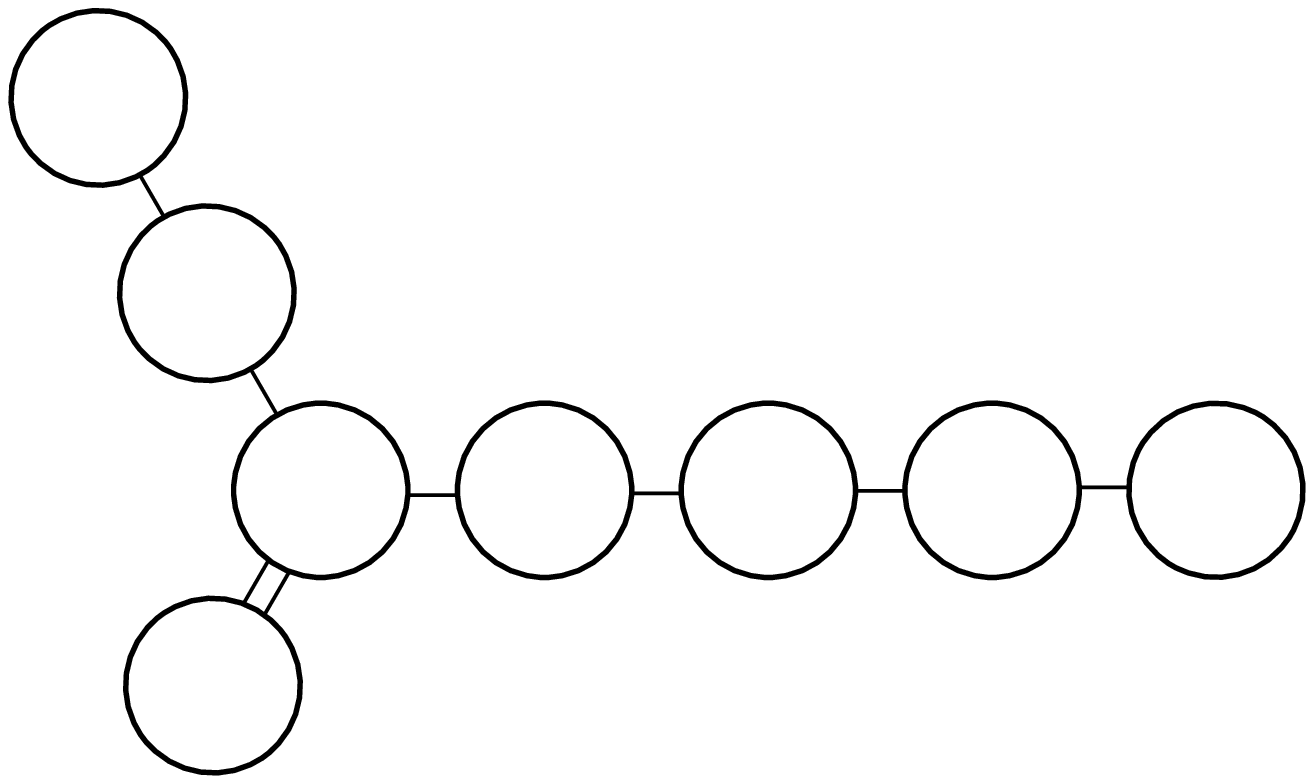}{.35} \kern 80pt
\raise10pt\hbox{\fig{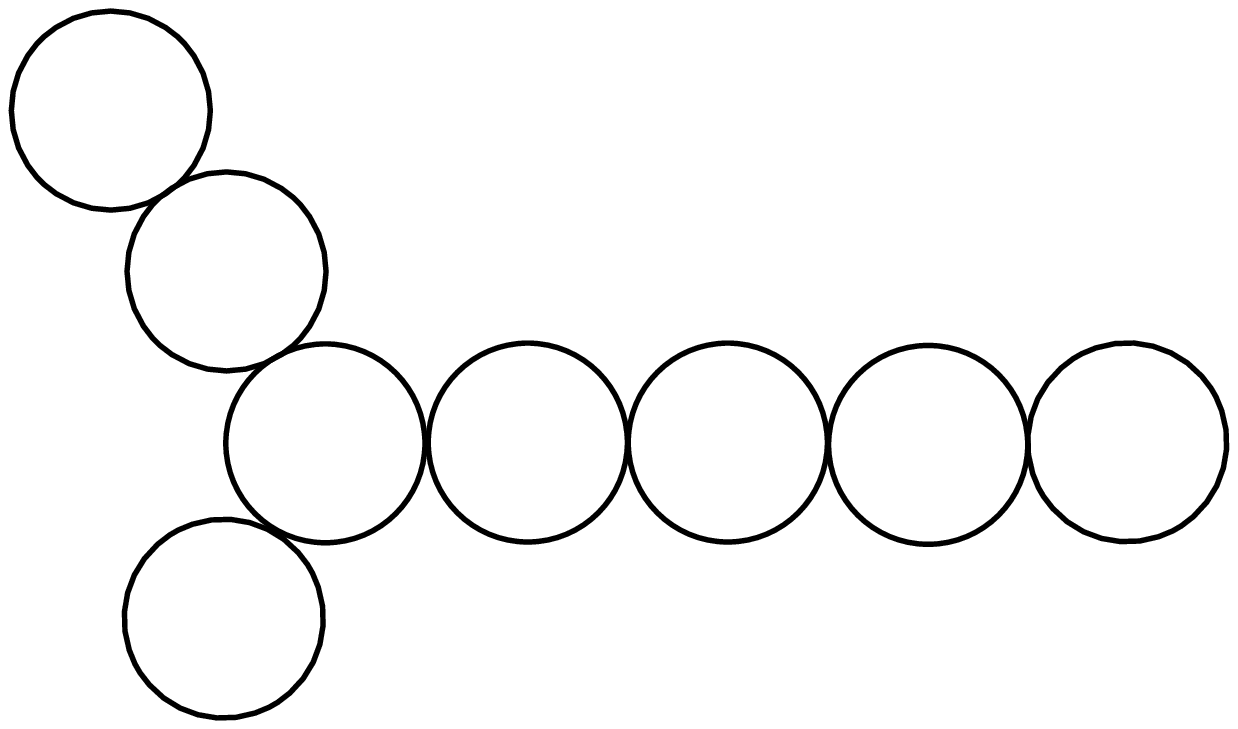}{.3}}} 
\smallskip
\centerline{\small\kern 50pt (a) \e--link with sticks for twists \kern 60 pt (b)
\e--link as grapes \kern 50pt}
\figno{1.1}
\end{figure}


\head{Grapes}

In Figure 1.1a the ``stick" notation \place{0}{-3}{\fig{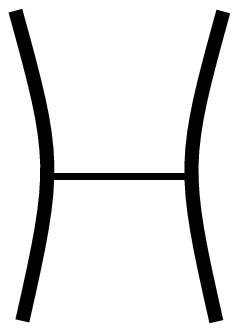}{.14}}\
\quad denotes a full left twist \place{0}{-3}{\fig{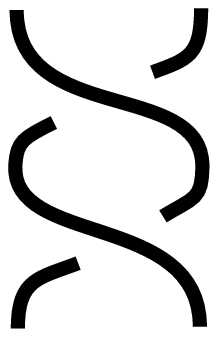}{.14}}\ \quad
between the vertical strands, while \place{0}{-3}{\fig{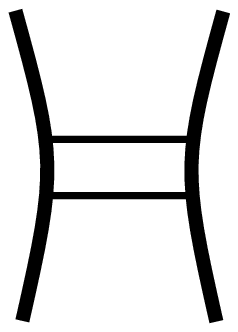}{.14}}\ \quad
denotes a full right twist \ \place{0}{-3}{\fig{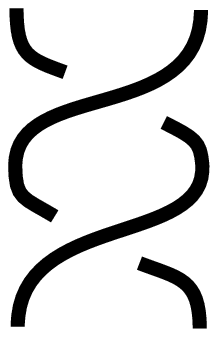}{.14}}\quad. 
The framings (when not labelled) should be assumed to be $-2$.  This
convention, which will be used throughout the paper, is convenient because
handle slides typically occur for 2--handles $A$ and $B$ with intersection
number $\pm1$ (where we identify a $2$--handle with its associated
homology class), so sliding $A$ over $B$ yields the class $A \pm B$ with
self-intersection $(A\pm B)^2 = A^2 \pm 2 A\cdot B + B^2 = -2+2-2 = -2$.

In  Figure 1.1b the circles are all required to lie in the hexagonal packing
of the plane. The convention is that a tangency between a pair of circles
represents a full twist between them, and that this twist is right-handed
if and only if the line joining their centers has positive slope.  (In
particular the drawings in Figure 1.1 represent identical link diagrams.)
Framings are $-2$ as in our standing convention.  Any configuration of
hexagonally packed circles representing a framed link (and thus a
$4$--manifold) by these conventions will be refered to as ``grapes"; each
individual circle will be called a grape.  It will be seen in what follows
that the framed links that arise in studying elliptic surfaces and related
$4$--manifolds can often be represented (after suitable handle-slides) by
grapes; this observation will streamline many of our constructions. 

Typical handleslides over a grape are illustrated in Figure 1.2, using the
stick notation.


\medskip
\centerline{\fig{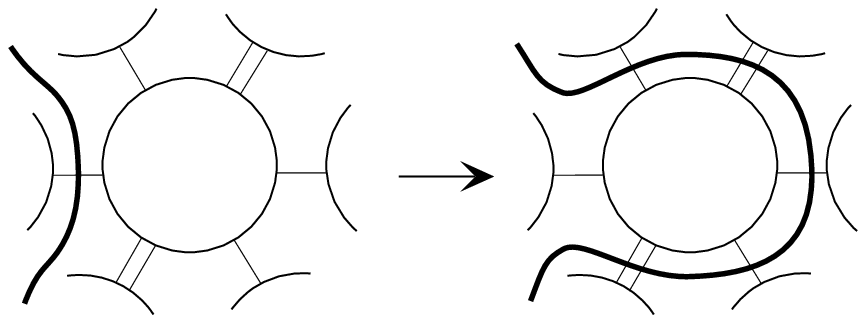}{.6} 
\kern 40pt \fig{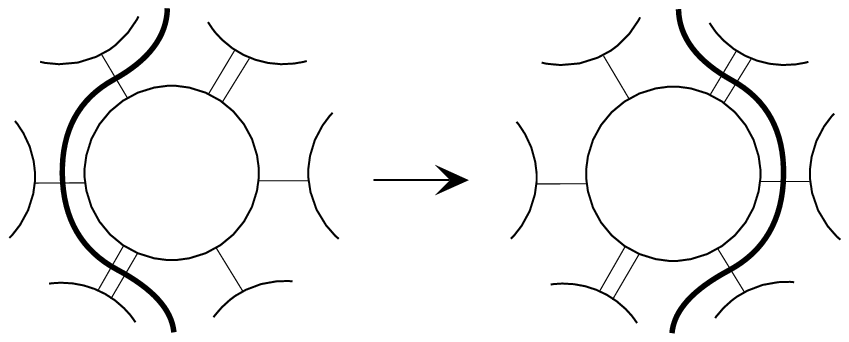}{.6}}
\centerline{\small(a) \kern 175 pt (b)}
\figno{1.2: Handleslides over a grape}

\medskip


\noindent
{\bf Branched covers}

\smallskip

The $4$--manifold $C_2$ which is the $2$--fold branched cover of $B^4$ along
a minimal genus Seifert surface $F_{3,5}$ for the $(3,5)$--torus
knot in $S^3$, pushed into the interior of $B^4$, is described in the next
figures.


\begin{figure}[ht!] 
\centerline{\fig{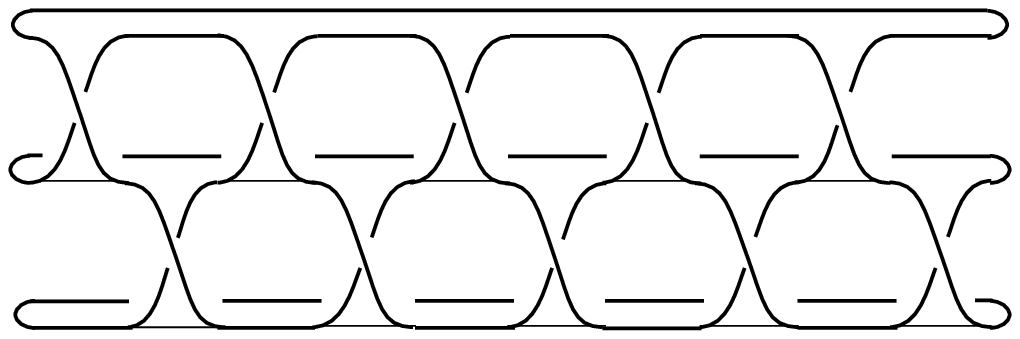}{.5} \kern 15pt 
\place{0}{25}{$=$}\kern 25pt \fig{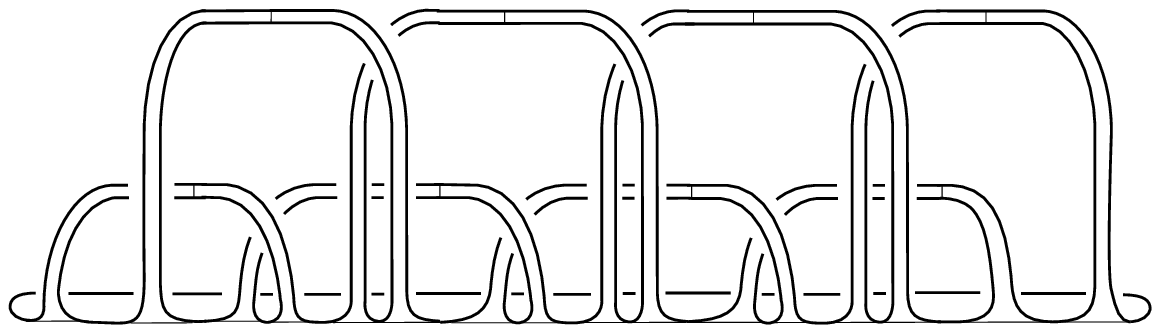}{.5}}
\bigskip
\centerline{\small(a) Seifert surface $F_{3,5}$}
\vskip .15in
\centerline{\fig{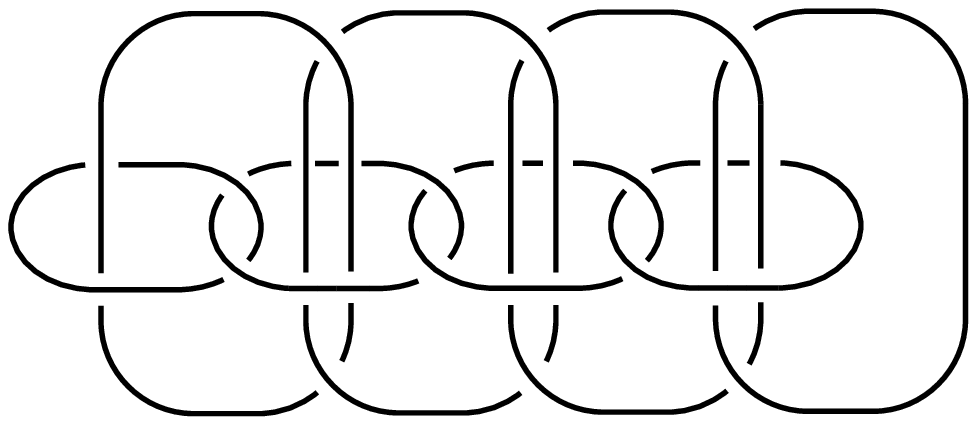}{.3} \kern 10pt
\place{3}{17}{$=$} \kern 25pt
\fig{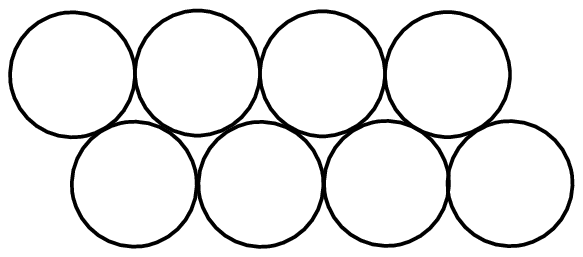}{.5} \kern 10pt}
\bigskip
\centerline{\small(b) Branched cover $C_2$ as a bunch of grapes}
\figno{1.3}
\end{figure}

Since torus knots are fibered, all Seifert surfaces of minimal
genus are isotopic, so the surface in Figure 1.3a will do.  Note that the Seifert surface in the first drawing consists of
three stacked disks with ten vertical half twisted bands joining them;  the
front four (large) $1$--handles in the second drawing come from the top disk
and the upper five half-twisted bands, the back four (small) $1$--handles
come from the middle disk and lower five half twisted bands, and the
$0$--handle is the bottom disk.

By the algorithm in \cite{AK} for drawing framed link descriptions of
branched covers of Seifert surfaces, a half circle should be drawn in each
$1$--handle, and then these eight half circles should be folded down to get
the link shown in Figure 1.3b.  The framing on each component is twice the
twist in corresponding $1$--handle (which is $-1$), so is $-2$ and not drawn
by convention.  Now folding the four smaller components over the top of the
larger ones gives the bunch of eight grapes shown in the second drawing.

In a similar way, we draw in Figure 1.4 the Seifert surface $F_{2,5}$ for
the $(2,5)$--torus knot, followed by its 3--fold branched cover $C_3$. 
In the algorithm (in \cite{AK}) {\sl two} half circles are drawn
in each of the four 1--handles, and then one set of four is folded down
followed by the other set.  This produces the first drawing in Figure
1.4b.  The second drawing is obtained from the first by sliding the outer
2--handle over the inner 2--handle for each of the four pairs of
2--handles.\foot{In fact there are two algorithms in \cite{AK} for drawing
the cover (see Figures 5 and 6 in \cite{AK}).  The first (which is more
natural but often harder to visualize) yields the grapes directly.  The
second is derived from the first by sliding handles --- the reverse of the
slides above in the present case.}  This link clearly coincides with the
bunch of grapes in Figure 1.3b, showing that $C_3$ is diffeomorphic to $C_2$.


\begin{figure}[ht!] 
\centerline{\fig{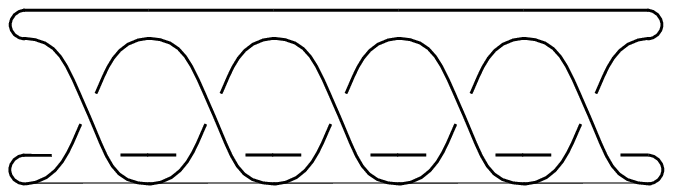}{.7} \kern 15pt \place{0}{15}{$=$} \kern 20pt
\fig{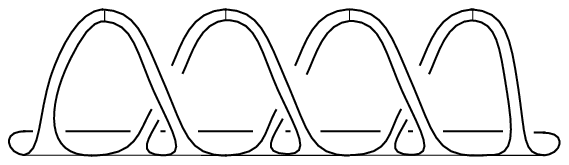}{.8}}
\bigskip
\centerline{\small(a) Seifert surface $F_{2,5}$}
\vskip .1 in 
\centerline{\fig{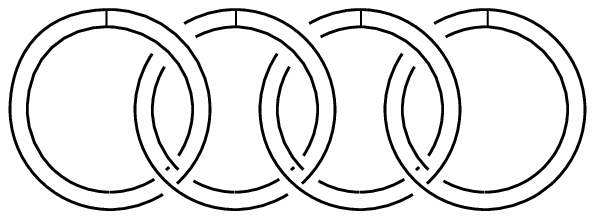}{.7} \kern 10pt \place{0}{15}{$\longrightarrow$}\kern
30pt \fig{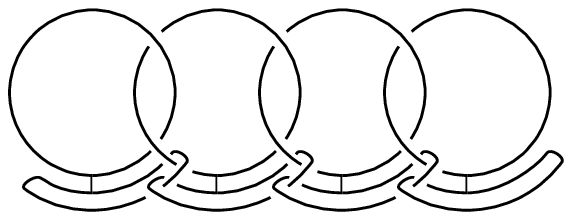}{.7}}
\bigskip
\centerline{\small(b) Branched cover $C_3$ (same bunch of grapes)}
\figno{1.4}
\end{figure}

Finally, Figure 1.5 shows the Seifert surface $F_{2,3}$ for the $(2,3)$--torus
knot and its 5--fold branched cover $C_5$
(where \place{0}{-3}{\fig{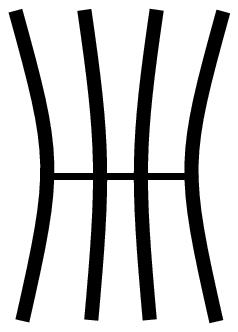}{.14}}\ \quad denotes a full left twist in
the vertical strands).  To pass from the first drawing in 1.5b to the second we
perform six handleslides, sliding each circle over its parallel neighbor,
starting with the outermost circles and working inward.  Rotating the last
drawing by a quarter turn yields the same bunch of grapes,$^{\dagger}$ showing
that  $C_5$ is diffeomorphic to $C_2$.


\begin{figure}[ht!] 
\centerline{\fig{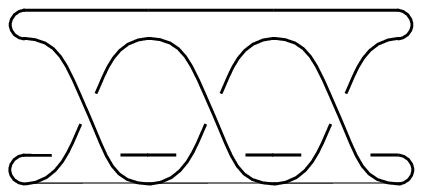}{.7} \kern 15pt \place{0}{15}{$=$} \kern 20pt
\fig{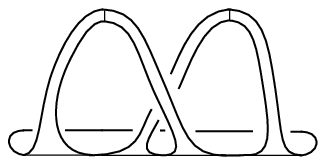}{.8}}
\bigskip
\centerline{\small(a) Seifert surface $F_{2,3}$}
\bigskip
\centerline{\raise15pt\hbox{\fig{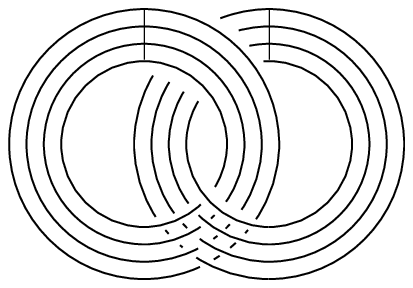}{.8}} \kern 20pt
\place{0}{40}{$\longrightarrow$}
\kern 40pt \fig{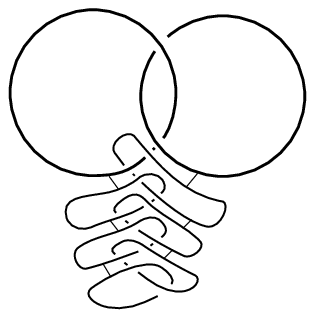}{.8} \place{20}{40}{$=$} \kern 50pt
\lower15pt\hbox{{\fig{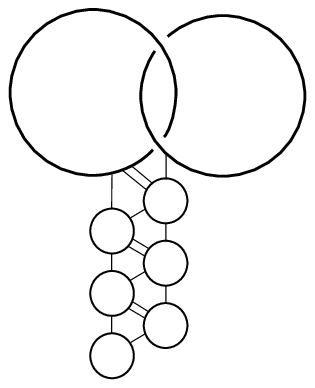}{.8}}}}
\centerline{\small(b) Branched cover $C_5$ (same bunch of grapes)}
\figno{1.5}
\end{figure}


\medskip\goodbreak
{\bf Equivalence of handlebody decompositions}
\medskip\nobreak

To show that these covers are diffeomorphic to \e, we introduce a move
on an arbitrary cluster of grapes (configuration of hexagonally packed
circles), called a {\it slip\/} \cite{LR}, which amounts to a sequence of
handle slides and isotopies: Suppose that such a cluster contains a grape
(labelled $A$ in Figure 1.6a) which is the first of a straight string of
grapes, in {\sl any of the six possible directions}.  If there are no grapes
in the dotted positions shown in Figure 1.5a, then grape $A$ can be moved by a
{\it slip\/} to the other end of the the string, that is to the position of
the dotted grape $B$ in the figure.  (Note that this slip can be reversed.)


\begin{figure}[ht!] 
\centerline{\fig{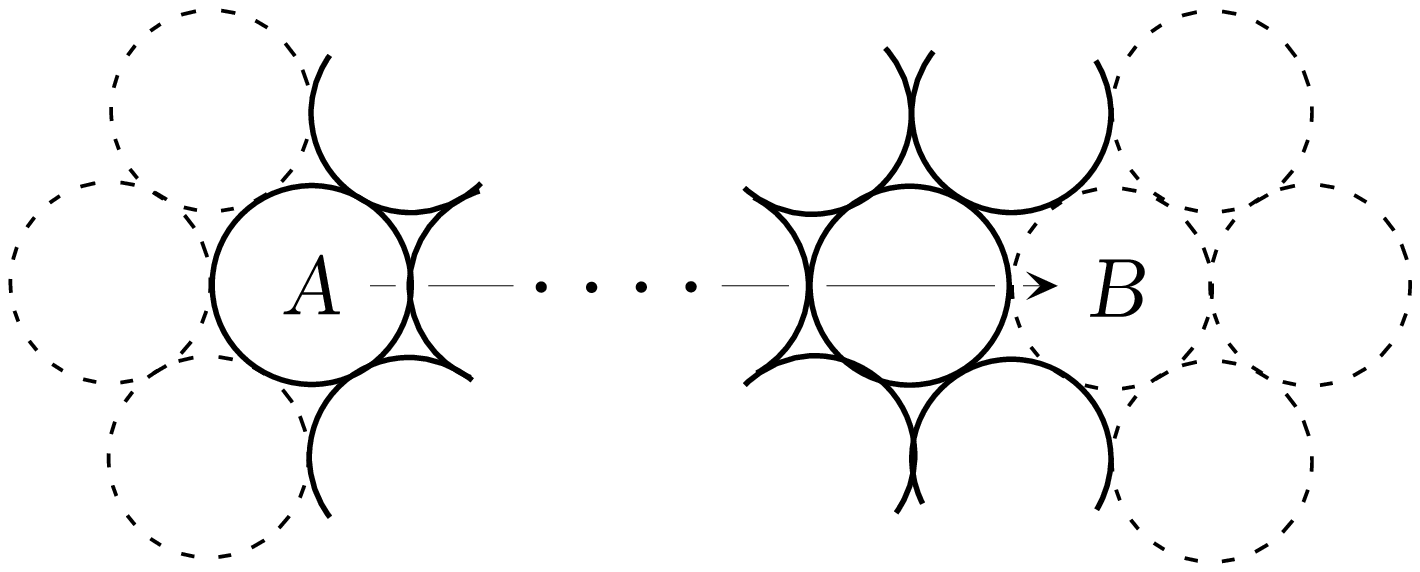}{.4}  \kern30pt  \place{-1}{0}{\fig{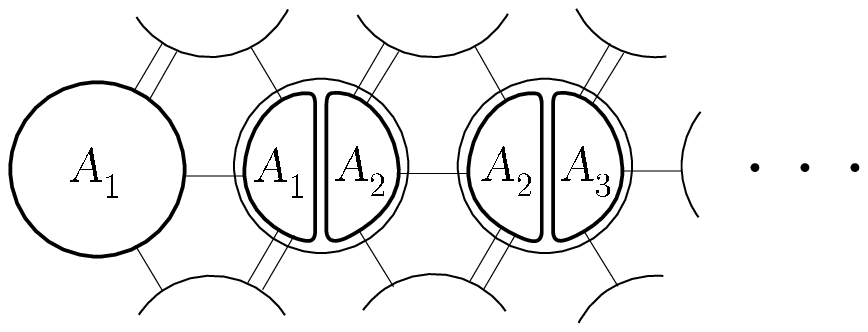}{.7}} 
\place{53}{26}{$^{'}$} \place{96}{26}{$^{'}$} \kern 170pt}
\bigskip
\centerline{\small\kern40pt(a) a slip \kern 120pt (b) the anatomy of a slip}
\figno{1.6}
\end{figure}

The handle slides and isotopies which produce the slip are indicated in
Figure 1.6b:\ $A=A_1 \rightsquigarrow A_1' \rightarrow A_2
\rightsquigarrow A_2' \rightarrow A_3 \rightsquigarrow \cdots
\rightarrow A_n \rightsquigarrow B$.  Here $A_i \rightsquigarrow A_i'$ is the
obvious isotopy (folding under when moving horizontally, folding
over when moving along the line inclined at $-\pi/3$ radians, and moving in
the plane of the paper by a regular isotopy when moving along the line
inclined at $\pi/3$ radians) as the reader can easily check, and
$A_{i+1}$ is obtained by sliding $A_i'$ over its encircling grape (cf\
Figure\ 1.2b).

Finally observe that the sequence of seven slips shown in Figure
1.7 takes the grapes defining \e\ (Figure\ 1.1b) to the grapes for the
branched cover $C_2$ (Figure\ 1.3b).
\bigskip


\centerline{\fig{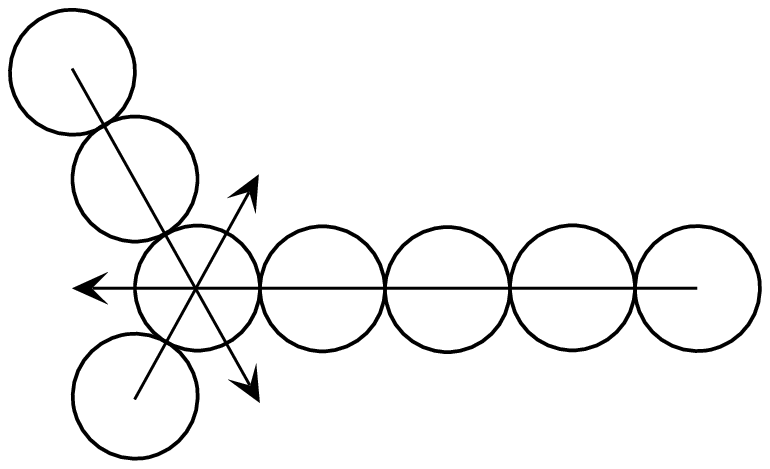}{.4} \place{10}{15}{$\longrightarrow$} \kern 40pt
\fig{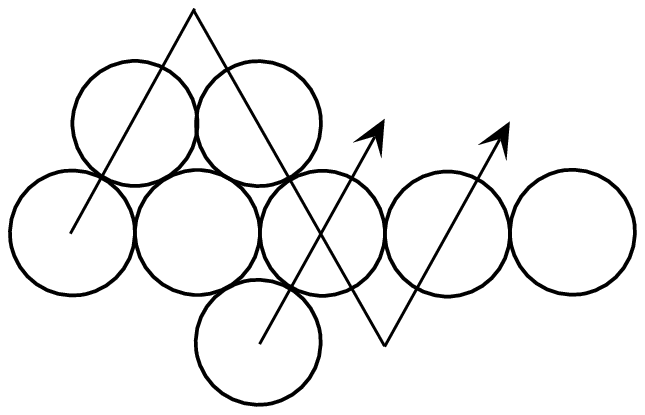}{.4}
\place{10}{15}{$\longrightarrow$} \kern 40pt \place{0}{5}{\fig{17c}{.4}}
\kern 60pt}
\figno{1.7: slippin' an' a slidin'}

Note that in the middle picture, the single slip must be
performed before the last leg of the triple slip.


\section{Singular fibers in elliptic surfaces}

The \e--manifold occurs naturally as a neighborhood of (most of) a singular
fiber in an elliptic surface, and so the discussion in the last section
suggests a general study of such neighborhoods.  We begin with a
brief introduction to the topology of elliptic surfaces and their singular
fibers; much fuller accounts for topologists appear in \cite{GS} and
\cite{HKK}.

An {\it elliptic surface\/} is a compact complex surface $E$ equipped with a
holomorphic map $\pi\co E\to B$ onto a complex curve $B$ such that the each
{\it regular fiber} (preimage of a regular value of $\pi$) is a
non-singular elliptic curve --- topologically a torus --- in $E$. Thus
$E$ is a $T^2$--bundle over $B$ away from the (finitely many) critical
values of $\pi$.  The fibers over these critical values are called the {\it
singular fibers\/} of the surface.  

Each singular fiber $C$ in $E$ is a union of irreducible curves
$C_i$, the {\it components\/} of the fiber. Topologically the components are
closed surfaces, possibly with self-intersection or higher order
singularities (for example a ``cusp"), and distinct components can
intersect, either transversely or to higher order.  Furthermore, each
component has a positive integer {\it multiplicity\/} $m_i$, where $\sum
m_iC_i$ represents the homology class of a regular fiber.  The multiplicity
of $C$ is then defined to be the greatest common divisor of the
multiplicities of its components.   We shall limit our discussion to {\it
simple\/} singular fibers, that is, fibers of multiplicity one.\foot{All
other fibers are either multiples of singular fibers or multiples of regular
fibers, and the latter have uninteresting neighborhoods, namely $T^2\times
B^2$.}  We also assume that $E$ is {\it minimal\/}, that is, not a blow-up of
another elliptic surface, or topologically not a connected sum of an
elliptic suface with $\overline{\C P^2}$.  This precludes any {\it
exceptional\/} components (non-singular rational curves of $-1$
self-intersection) in singular fibers.   

The singular fibers in minimal elliptic surfaces were classified by
Kodaira \cite{Kod}, and the simple ones fall into eight classes: two infinite
families I$_n$ and I$^*_n$ (where $n$ is a non-negative integer, positive in
the first case since I$_0$ represents a regular fiber), three additional
types II--IV and their ``duals" II*--IV* (explained below). In all cases the
components are rational curves --- topologically 2--spheres --- and so the
singular fiber can be depicted by a {\sl graph} of intersecting arcs
representing these components.


\head{Fibers of type \  I--IV \ (Table 1)}

The two simplest singular fibers are the {\it fishtails} (type I$_1$)
and the {\it cusps} (type II).  A fishtail consists of a single
component, an immersed 2--sphere with one positive double point, and is
represented by a self-intersecting arc \,\fig{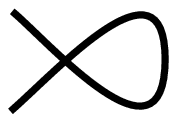}{.2}\,.  A cusp is a
2--sphere with one singular point which is locally a cone on a right handed
trefoil knot; this is denoted by a cusped arc
\lower1pt\hbox{\,\fig{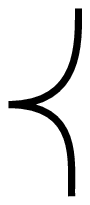}{.2}\,}.  

The components in all other singular fibers (including those of $*$--type) are
smoothly embedded 2--spheres with self-intersection $-2$.  In particular, a
singular fiber of type III consists of two such 2--spheres which are tangent
to first order at one point, denoted by a pair of tangent arcs
\,\fig{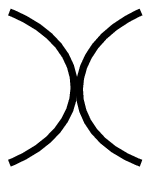}{.2}\,, and one of type IV consists of three such 2--spheres
intersecting transversely in one triple point, denoted \,\fig{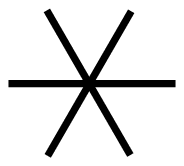}{.2}\,. 
Singular fibers of type I$_n$ for $n>1$, called {\it necklace fibers},
consist of $n$ such 2--spheres arranged in a cycle, each intersecting the one
before it and the one after it (which coincide if $n=2$).  For example
I$_2=$ \,\fig{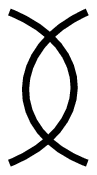}{.2}\, and I$_5=$ \,\fig{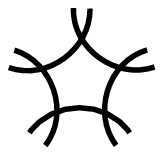}{.25}\,.  These graphs
are reproduced in the first column of Table 1 below; note that all
components have multiplicity one in these types of fibers.

The second column in the table gives natural framed link descriptions for
regular neighborhoods of these fibers, following \cite{HKK}.  A
neighborhood of a fishtail is clearly a self-plumbing of the cotangent
disk bundle $\tau^*$ of the $2$--sphere, of euler class $-2$ (the
homology class represented by the fishtail must have self-intersection
zero, so the euler class is $-2$ to balance the two positive points of
intersection arising from the double point) and this can be constructed
as a $0$--handle with a round $1$--handle (\,= 1--handle plus a 2--handle)
attached as shown.  Similarly a neighborhood of an $n$--component necklace
fiber is a circular plumbing of $n$ copies of $\tau^*$, or equivalently
surgery on a chain of circles in $S^1\times B^3$ (\,= 0--handle plus a
1--handle).  The cusp neighborhood is obtained by attaching a single 2--handle
along the zero framed right-handed trefoil. Fibers of type II and III are
gotten by attaching handles to the $(2,4)$ and $(3,3)$--torus links,
respectively, with  $-2$ framings on all components.  (In \S3 we will give
explicit models for these neighborhoods in which the projection of the
elliptic surface is evident.)


\begin{figure}[ht!]\small
\centerline{\begin{tabular}{c||c|c|c|} 
type & graph &\vrule height 2pt depth 8pt width 0pt framed link & monodromy \\ 
\hline \hline 
&&& \\
I$_1$ & \lower5pt\hbox{\fig{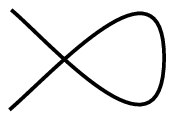}{.4}} & \lower8pt\hbox{\fig{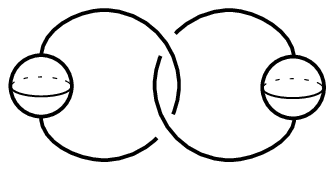}{.45}
\place{-16}{16}{$^0$}} & $V 
$ \\  
&&& \\
\hline 
&&& \\
I$_n \ (n\ge2)$ & \lower5pt\hbox{\fig{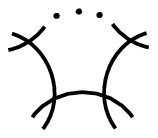}{.4}} &
\ \ \lower8pt\hbox{\fig{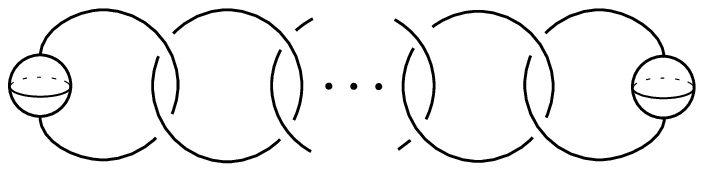}{.45} \place{-72}{17}{$^{-2}$}
\place{-42}{17}{$^{-2}$} \place{-25}{17}{$^{-2}$}} & $V^n 
$ \\  
&&& \\ 
\hline 
&&& \\
II & \lower5pt\hbox{\fig{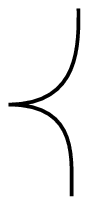}{.4}} & \lower8pt\hbox{\fig{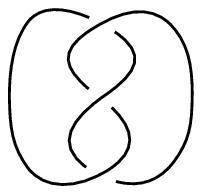}{.5}
\place{-2}{12}{$^0$}} & $UV 
$ \\ 
&&& \\ 
\hline 
&&& \\
III & \lower5pt\hbox{\fig{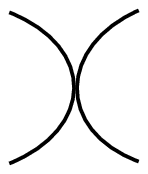}{.4}} & \lower8pt\hbox{\fig{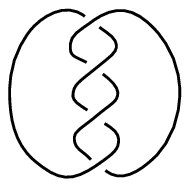}{.5}
\place{-42}{11}{$^{-2}$} \place{-6}{11}{$^{-2}$}} &
$UVU 
$ \\  
&&& \\ 
\hline 
&&& \\
IV & \lower5pt\hbox{\fig{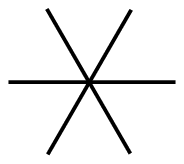}{.4}} & \lower8pt\hbox{\fig{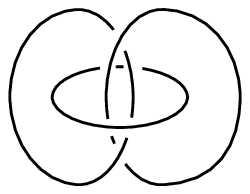}{.5}
\place{-49}{12}{$^{-2}$} \place{-20}{0}{$^{-2}$}
\place{-8}{12}{$^{-2}$}} &
$(UV)^2 
$ \\  &&& \\
\hline
\end{tabular}}
\medskip
\tabno{1}{Singular fibers of type I--IV}
\end{figure}

The final column in the table gives the {\it monodromy\/} of the torus
bundle around each singular fiber with respect to a suitably chosen
basis for the first homology of a regular fiber,\foot{More precisely,
if we pick a base point $b_0$ in $B$ and a basis for the first homology
of the fiber over $b_0$, choose paths connecting $b_0$ to each critical
value $p_i$ of $\pi$, and choose small loops $\gamma_i$ around each
$p_i$, then we get a well defined $2\times 2$--matrix for each singular
fiber, representing the monodromy of the torus bundle over the
associated $\gamma_i$.} given in terms of the generators
$$
U = \pmatrix \ 1&0\\-1&1 \endpmatrix \qquad{\rm and}\qquad 
V = \pmatrix 1&1\\0&1 \endpmatrix 
$$
of $SL(2,\Z)$.  Note that $(UV)^6 = I = (UVU)^4$ (since $UVU=VUV$).  Also
note that if the chosen basis is viewed as a {\it longitude\/} and {\it
meridian\/} of the regular fiber (in that order) then $U$ corresponds to
meridianal Dehn twist, and $V$ to a longitudinal one.

The fishtail neighborhood can also be obtained from a thickened regular
fiber $N = T^2\times B^2$ by attaching a $2$--handle with framing $-1$ to an
essential embedded circle $C$ (or {\it vanishing cycle\/}) lying in a torus
fiber in $\partial N = T^2\times S^1$ (see for example \cite{Kas}).  This
changes the trivial monodromy of $\partial N$ by a Dehn twist about $C$,
giving $V$ for the monodromy of the fishtail if $C$ is the longitude in the
torus.  Figure 2.1a shows the standard handlebody decomposition of $N$
with two 1--handles and a 0--framed ``toral" 2--handle (where for convenience
we identify the horizontal and vertical directions with the meridian and
longitude on the torus fiber), and Figure 2.1b shows the result after
attaching the last 2--handle along a vertical (longitudinal) vanishing
cycle.  This handlebody will be denoted by $N_V$, and simplifies to the
one in the table by cancelling the vertical 1 and 2--handles.


\begin{figure}[ht!]
\centerline{\fig{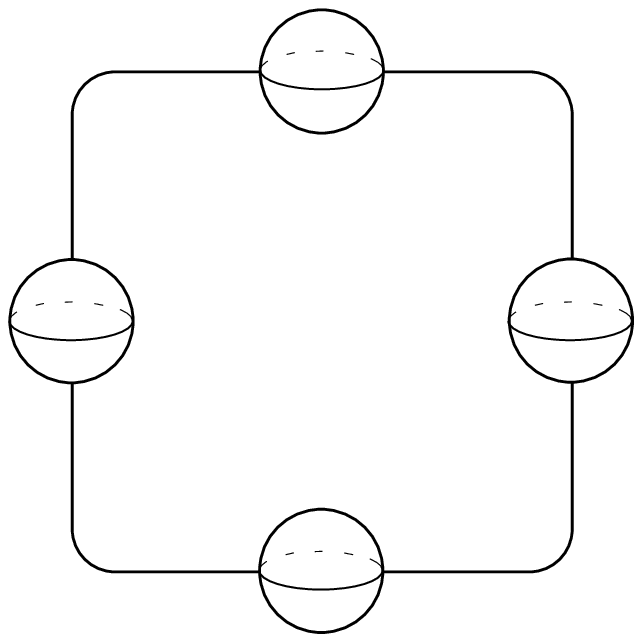}{.5} \place{-92}{70}{$^0$} 
\kern 50pt \fig{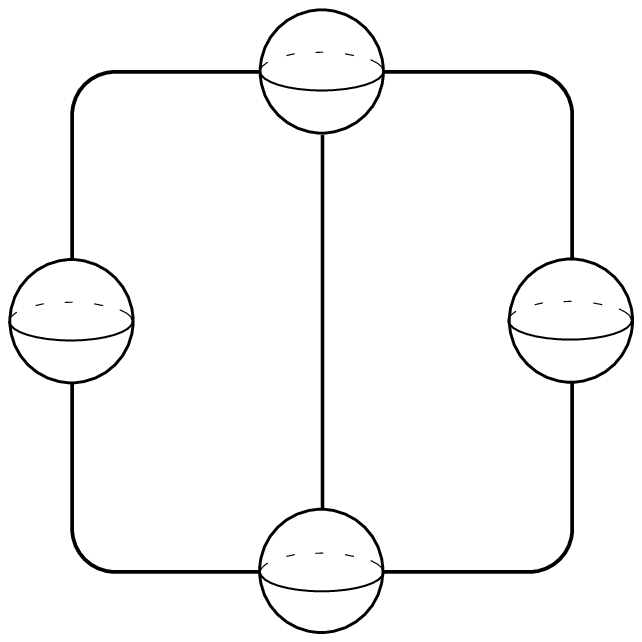}{.5}
\place{-92}{70}{$^0$}\place{-60}{45}{$^{-1}$}}
\medskip
\centerline{\small(a) $N = T^2\times B^2$ \kern 60pt (b) $N_V =$ Fishtail}
\figno{2.1}
\end{figure}

Now it is well known that any simple singular fiber in an elliptic surface
breaks up into finitely many fishtails under a generic perturbation of the
projection near the fiber.  To show this explicitly for the fibers of type
I--IV (the argument for the other types will be given later) observe that
the factorization of the monodromy given in Table 1 suggests a pattern of
vanishing cycles.  For the necklace fiber I$_n$ with monodromy $V^n$ one
expects $n$ longitudes.  For the remaining types with monodromies
$UV\cdots$ the vanishing cycles should alternate between meridians and
longitudes.

More precisely, for any word $W$ in $U$ and $V$, consider the handlebody
$N_W$ obtained from $N$ by attaching a sequence of 2--handles along
$-1$--framed meridians (for each $U$ in $W$) and longitudes (for each
$V$) in successive torus fibers in $\partial N$.  Then we have the
following result (cf\ Theorem 1.25 in \cite{HKK}).

\begin{theorem} A regular neighborhood of a singular fiber of type {\rm
I$_n$, II, III} or {\rm IV} is diffeomorphic to $N_{V^n}$, $N_{UV}$,
$N_{UVU}$ or $N_{(UV)^2}$, respectively. 
\end{theorem}

Before giving the proof, we describe a general procedure for simplifying
$N_W$, illustrated with the word $W=(UV)^3V$ (which arises as the
monodromy of a fiber of type I$^*_1$ below).  The associated handlebody is
shown in Figure 2.2a.  Sliding each vanishing cycle over its parallel
neighbor, working from the bottom up, produces a bunch of grapes (all with
framings $-2$ as usual) hanging from the top horizontal and vertical cycles
(Figure\ 2.2b) --- the same process was used to identify the branched covers
$C_2$ and $C_5$ in the last section (Figure\ 1.5b).  Now cancelling the
2--handles attached to these last two cycles with the 1--handles gives the
handlebody $R_W$ in Figure 2.2c.\foot{Note that if one labels the upper row
of grapes with $U$ and the lower row with $V$, as shown, then reading from
left to right yields the truncation of $W$ obtained by deleting the initial
$UV$.}  This process will be called the {\it standard reduction\/} of the
handlebody $N_W$ to the {\it reduced form} $R_W$.


\begin{figure}[ht!] 
\centerline{\kern 20pt \fig{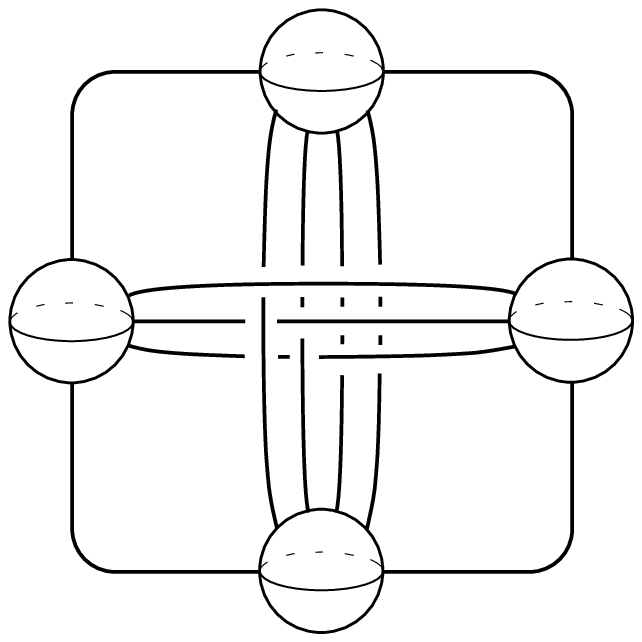}{.5} \place{-92}{70}{$^0$}
\place{-75}{60}{$\scriptstyle\textup{all}$}
\place{-82}{53}{$\scriptstyle -1\textup{'s}$} 
\kern 60pt \fig{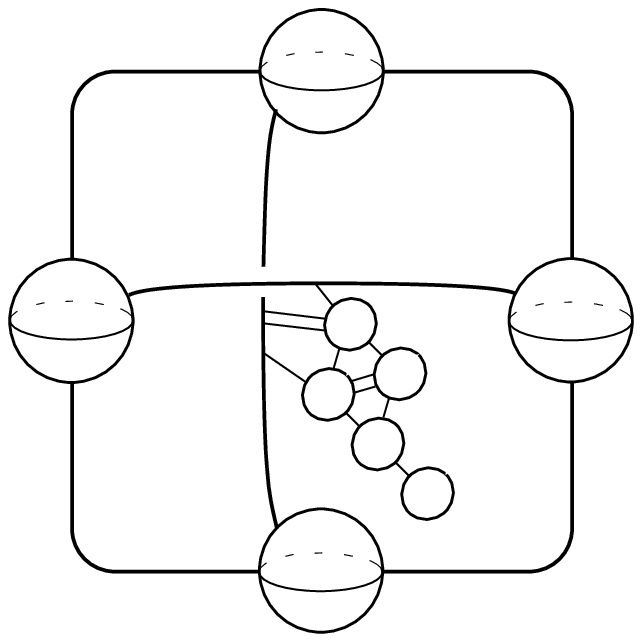}{.5}
\place{-92}{70}{$^0$}\place{-68}{60}{$^{-1}$}\place{-40}{48}{$^{-1}$}}
\medskip
\centerline{\small(a) The handlebody $N_{(UV)^3V}$ \kern 40pt (b) Slide to grapes}
\bigskip
\centerline{\kern 50pt \fig{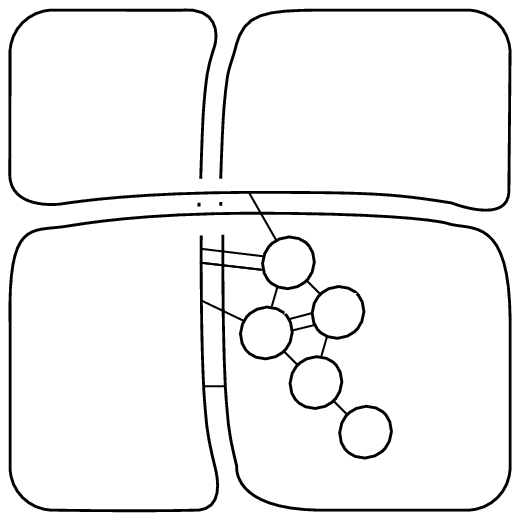}{.5} \place{-87}{55}{$^0$} 
\place {30}{30}{$\rightsquigarrow$} \place {18}{20}{{\small isotopy}} \kern
75pt \fig{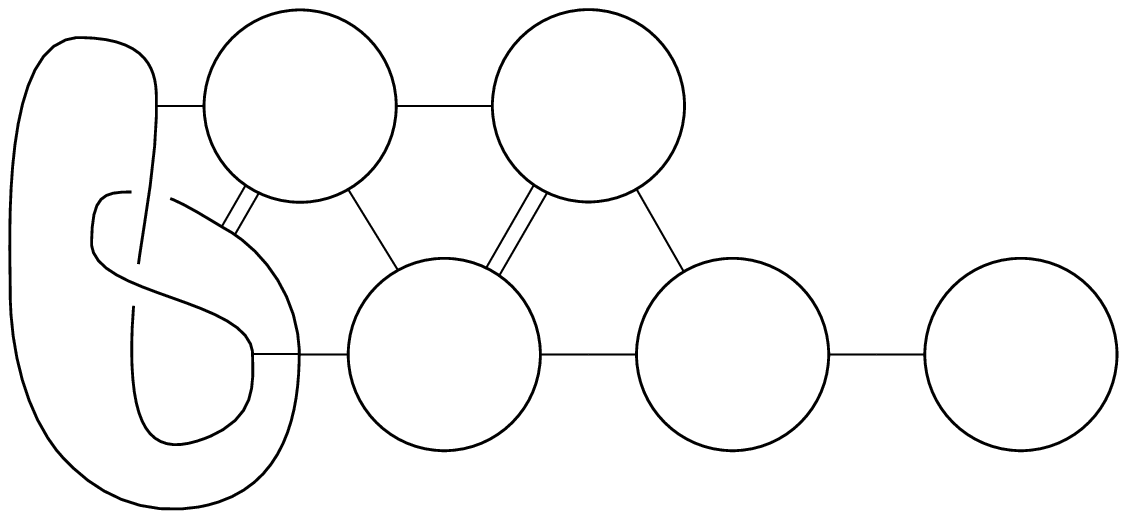}{.5}
\place{-199}{50}{$^0$}
\place{-156}{52}{$^U$} \place{-117}{52}{$^U$} 
\place{-141}{15}{$^V$} \place{-102}{15}{$^V$} \place{-62}{15}{$^V$}}
\medskip
\centerline{\small(c) Cancel 1--handles to get reduced form $R_{(UV)^3V}$}
\figno{2.2: Standard reduction}
\end{figure}

Observe that the first step of the reduction of $N_W$ (producing the grapes)
can be carried out in general, but that the cancellation of {\sl both}
1--handles requires at least one $U$ and one $V$ in $W$ (and the final
picture will look slightly different if $W$ does not start with $UV$).  In
particular for $W = V^n$ the reduced form $R_{V^n}$ is obtained by cancelling
the vertical 1--handle with the last vanishing cycle, as shown in Figure 2.3.


\begin{figure}[ht!]
\centerline{\kern 30pt \fig{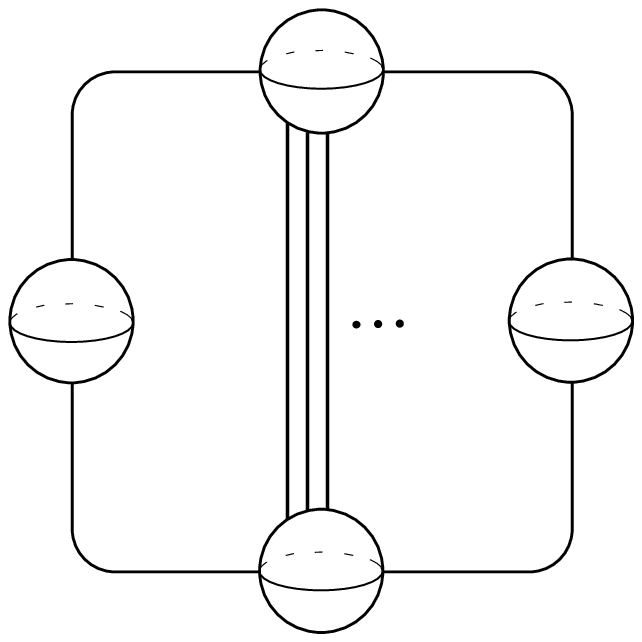}{.5} \place{-92}{70}{$^0$}
\place{-70}{60}{$\scriptstyle\textup{all}$}
\place{-77}{53}{$\scriptstyle -1\textup{'s}$}
\kern 60pt \raise18pt\hbox{\fig{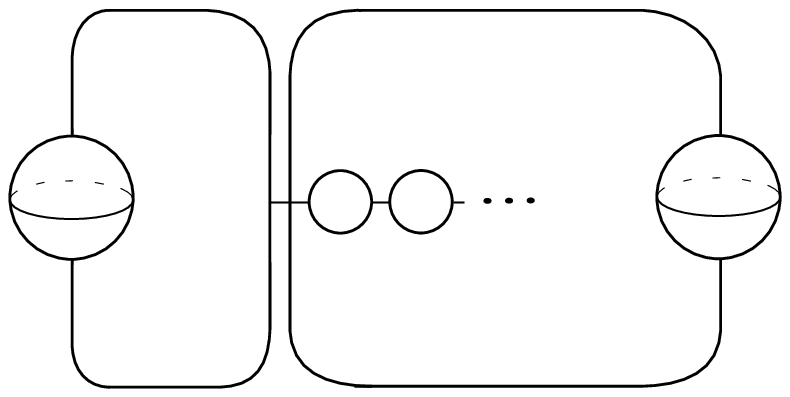}{.5} \kern10pt}
\place{-127}{60}{$^0$}}
\smallskip
\centerline{\small(a) The handlebody $N_{V^n}$ \kern 60pt (b) Reduced form
$R_{V^n}$}
\figno{2.3}
\end{figure}

We now return to the proof of the theorem.  For the necklace fiber I$_n$, the
obvious handleslides of the right-hand loop of the toral 2--handle in
$R_{V^n}$ (Figure\ 2.3b) over the grapes yields the picture for the
neighborhood of I$_n$ given in Table 1.  For the cusp, the reduced form
$R_{UV}$ (shown in Figure 2.4a) is exactly the 0--framed trefoil.  For fibers
of type III and IV, handleslides of the toral 2--handle in $R_{UVU}$ and
$R_{(UV)^2}$ are indicated in Figure 2.4, and an isotopy in each case yields
the corresponding picture in the table.


\begin{figure}[ht!] 
\centerline{\kern30pt \fig{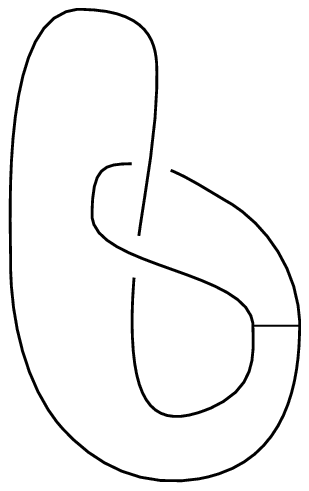}{.5} \place{-50}{50}{$^0$} \kern 60pt
\fig{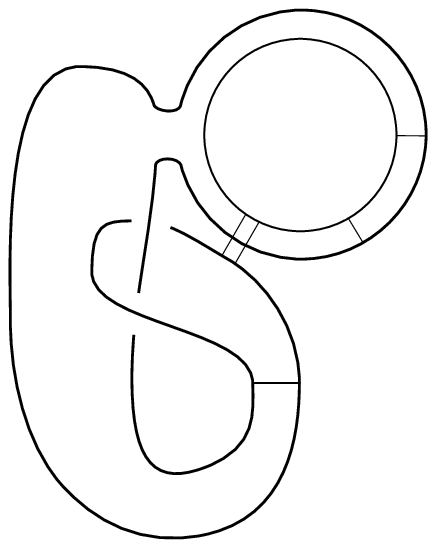}{.5} \kern 40pt \fig{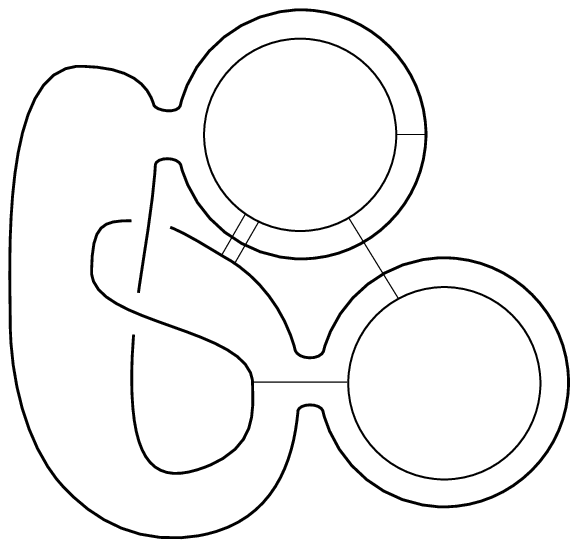}{.5}}
\medskip
\centerline{\small(a) II \kern 80pt (b) III \kern 80pt (c) IV}
\figno{2.4}
\end{figure}

\noindent This completes the proof of Theorem 1.
\qed


\head{Fibers of type \  I*--IV* \ (Table 1*)}

For the  singular fibers of $*$--type, any pair of components intersect
transversely in at most one point and there
are no ``cycles" of components.  Thus it is customary to represent
these fibers by the {\sl dual tree} with a vertex for each component $C_i$
and an edge joining any two vertices whose associated components intersect. 
The multiplicities $m_i$ of the components, which are often greater than one,
are recorded as weights on the vertices of the tree.  Note that the $m_i$ are
uniquely determined by the equations $\sum_i m_iC_i \cdot C_j = 0$
since $\sum_i m_i C_i$ is homologous to a regular fiber which is disjoint from
$C_j$; this translates into the condition that the weight of each vertex is
half the sum of the weights of its neighboring vertices.  

Regular neighborhoods of these fibers are the associated plumbings of the
cotangent disk bundle of $S^2$ (see \S3 for explicit models which explain
why these plumbings appear) and can all (with the exception of I$^*_0$) be
represented by grapes as in \S1, in fact in a variety of ways.  We choose
one, and record it in Table 1* along with the associated weighted tree and
monodromy.  

\begin{figure}[ht!] \small
\centerline{\begin{tabular}{c||c|c|c|} 
type & weighted tree &\vrule height 2pt depth 8pt width 0pt grapes 
& monodromy \\ 
\hline \hline 
I$^*_0$ & \zero & \zerog & $(UV)^3 = -I$ \\  
&&& \\
\hline 
I$^*_n \ (n>0)$ & \bee & \beeg & $(UV)^3V^n = -V^n$ \\  
&&& \\
\hline 
II* & \ \two \ & \twog & $(UV)^5 = (UV)^{-1}$ \\ 
&&& \\
\hline 
III* & \three & \threeg & $(UV)^4U = (UVU)^{-1}$ \\  
&&& \\
\hline 
IV* & \four & \fourg &$(UV)^4 = (UVUV)^{-1}$ \\  
&&& \\
\hline
\end{tabular}}
\medskip
\tabno{1*}{\small Singular fibers of type I*--IV*}
\end{figure}

The reader should note the inverse relation between the monodromies of each
fiber of type II--IV and its starred counterpart, which provides one
explanation for their common label.  This implies that neighborhoods of dual
fibers (ie,\ II and II*, III and III*, or IV and IV*) can be
identified along their boundaries to form a (closed) elliptic surface. 
Indeed one obtains in this way various non-generic projections of the
rational elliptic surface $E(1)$, diffeomorphic to $\C P^2\#9\overline{\C
P^2}$, whose generic projection has twelve fishtails (see eg\
\cite[\S1]{HKK}).

\nopagebreak
We now state the analogue of Theorem 1.

\setcounter{theorem}{0}

\begin{theorem}{${\kern-8pt^*}$} A regular neighborhood of a singular fiber
of type {\rm I$^*_n$, II*, III*} or {\rm IV*} is diffeomorphic to
$N_{(UV)^3V^n}$, $N_{(UV)^5}$, $N_{(UV)^4U}$ or $N_{(UV)^4}$, respectively. 
\end{theorem}

\begin{proof} Proceeding as in the proof of Theorem 1, consider the
handleslides of the toral 2--handle over the grapes in the reduced forms, as
shown in Figure 2.5.


\begin{figure}[ht!]
\centerline{\kern 5pt \fig{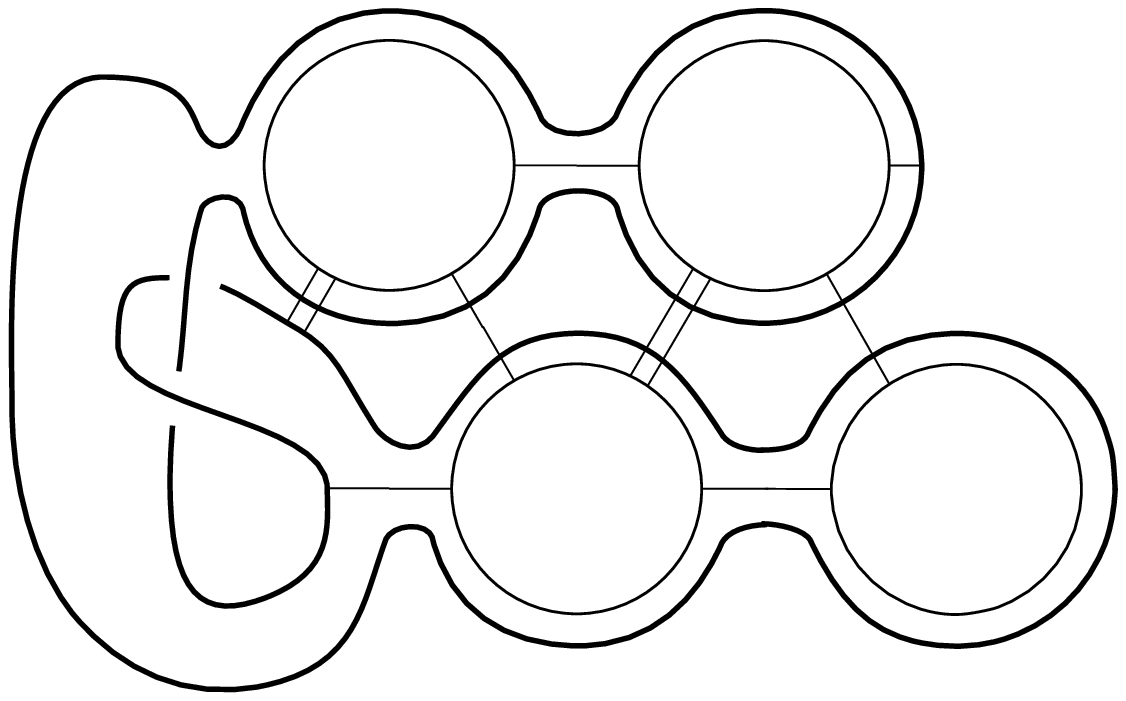}{.35} \kern 25pt \fig{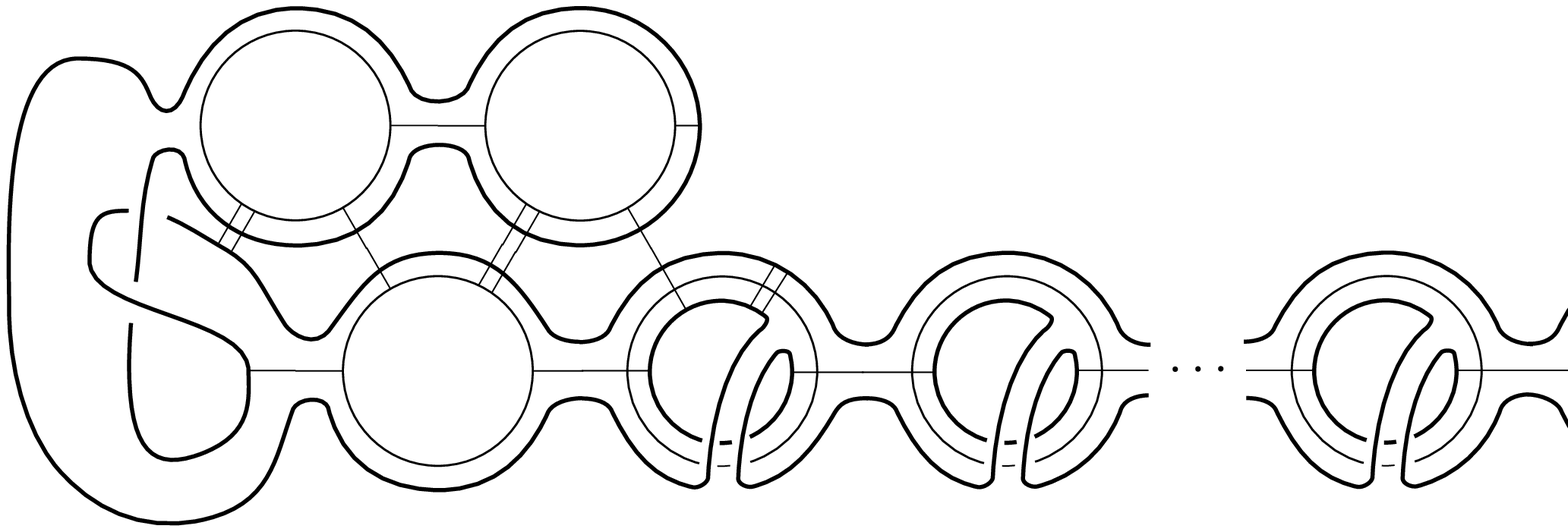}{.35}}
\medskip
\centerline{\small(a) I$^*_0$ \kern 150pt (b) I$^*_n$, $n>0$ \kern 40pt}
\bigskip
\centerline{\fig{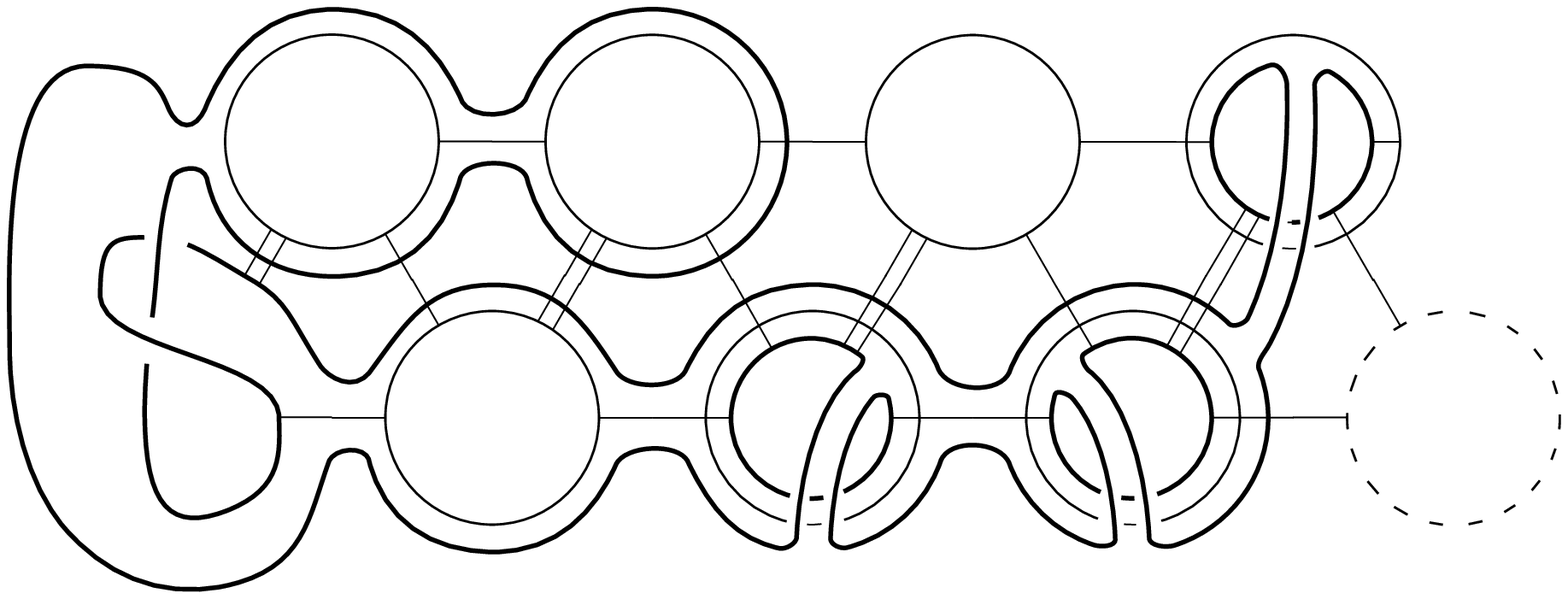}{.35} \kern 20pt \fig{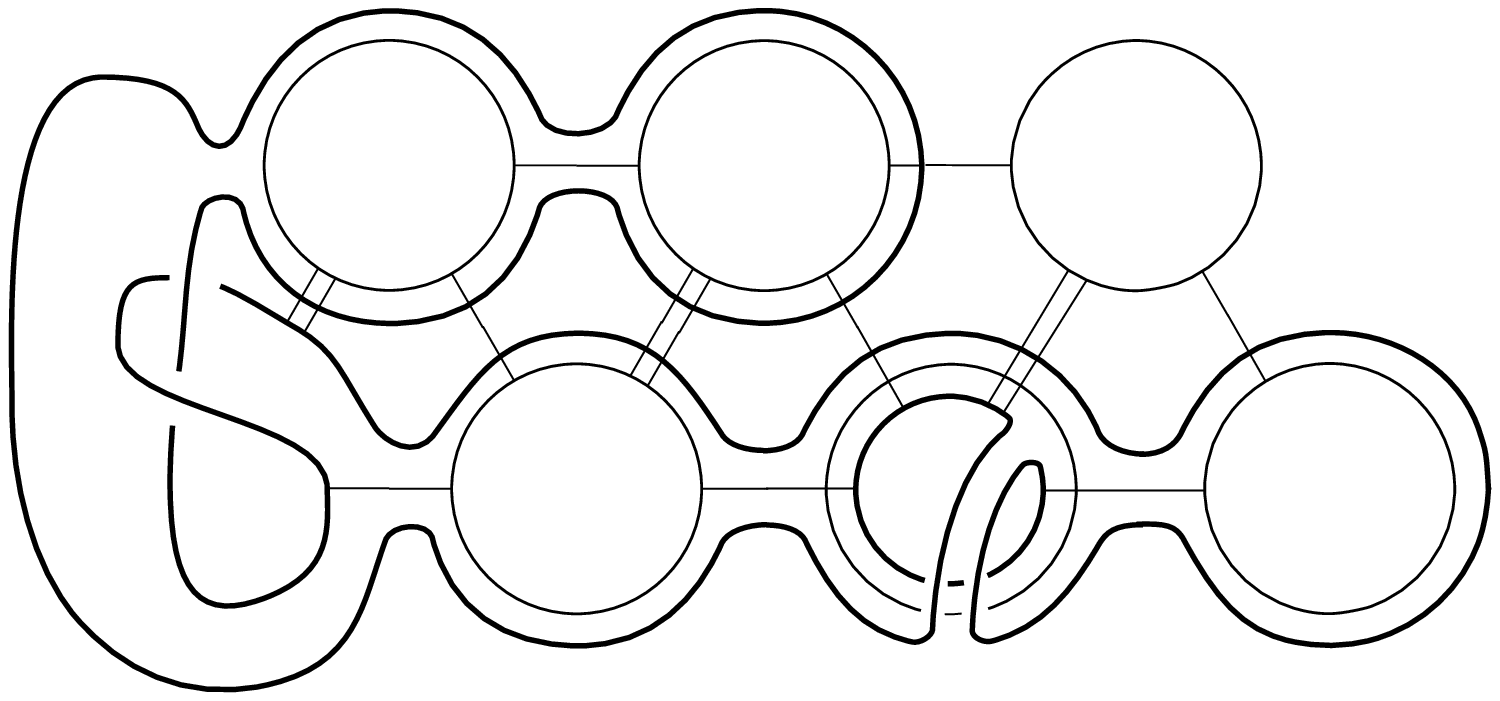}{.35}}
\medskip
\centerline{\small(c) II* or III* \kern 140pt (d) IV*}
\figno{2.5}
\end{figure}

After an isotopy, the toral 2--handle appears as in Figure 2.6,
labelled with a $T$; note that its framing is now
$-2$.  With the exception of I$_0^*$, the grapes are written in the usual
hexagonally packed notation.


\begin{figure}[ht!]
\centerline{\raise5pt\hbox{\fig{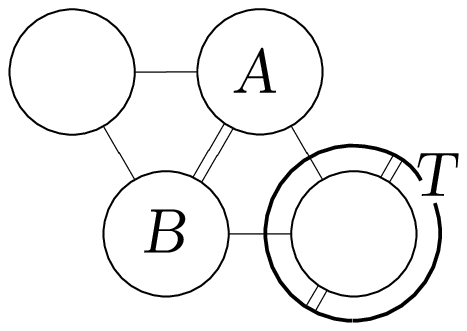}{.5}} \kern 60pt
\fig{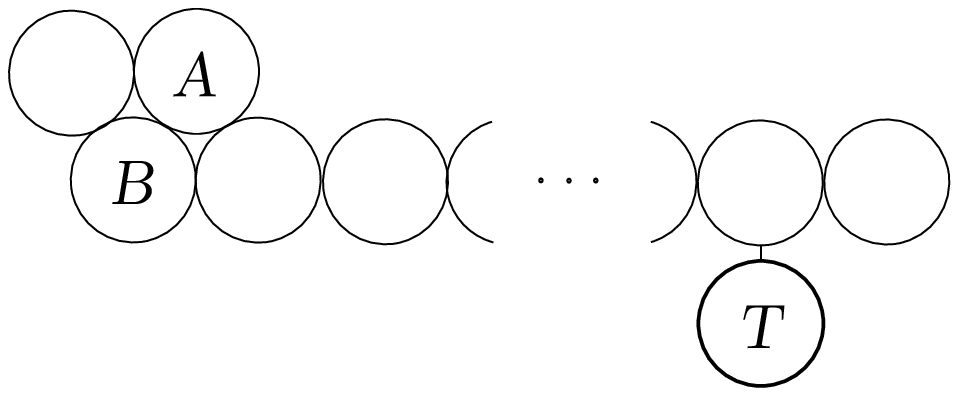}{.5}}
\smallskip
\centerline{\small(a) I$^*_0$ \kern 110pt (b) I$^*_n$, $n>0$ \kern20pt}
\bigskip\medskip
\centerline{\raise5pt\hbox{\fig{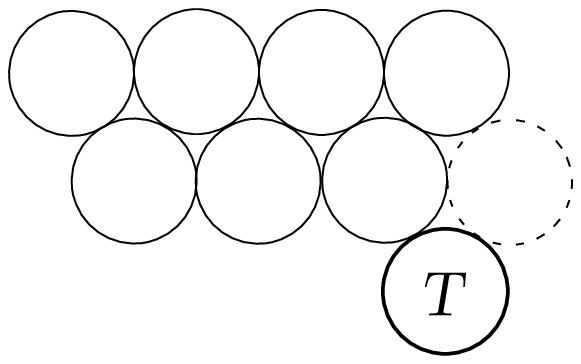}{.5}} \kern 80pt
\fig{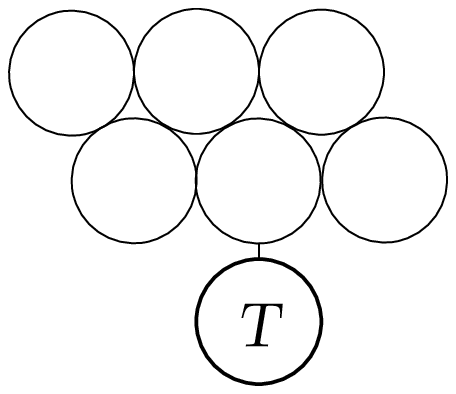}{.5}}
\smallskip
\centerline{\small(c) II* or III* \kern 100pt (d) IV*}
\figno{2.6}
\end{figure}

We should remark that one is guided in discovering the slides in Figure 2.5
by the multiplicities of the components of the singular fibers, which
determine the number of times the toral handle must slide over each grape to
achieve the simple pattern in Figure 2.6.  However, one must still be very
careful in choosing where to perform the slide in order to avoid
knotting and linking.

The argument is now completed by a sequence of slips.  For singular fibers
of type I* (Figure\ 2.6a,b) a single slip of grape $A$ over grape $B$ does
the job, recovering the handlebodies given in Table 1*.  For types II*--IV*,
the sequence of slips shown in Figure 2.7 will give the clusters of grapes given
in the table.  Note that if $T$ links only one of the grapes, then that grape
can be slipped over others while carrying $T$ along for the ride.


\begin{figure}[ht!] 
\centerline{\fig{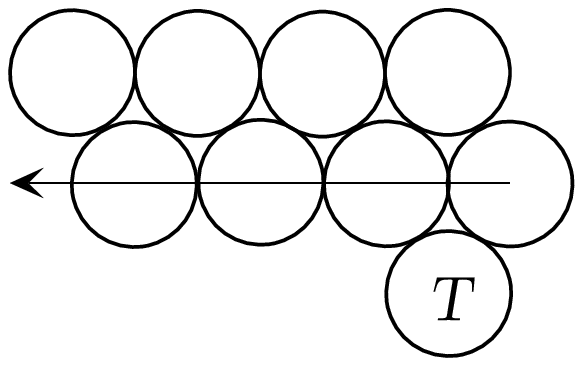}{.35}
\place{2}{15}{$\longrightarrow$}
\kern 25pt
\fig{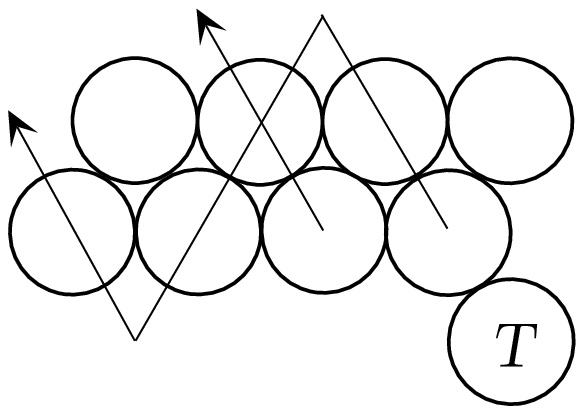}{.35} 
\place{7}{15}{$\longrightarrow$} 
\kern 25pt
\lower5pt\hbox{\fig{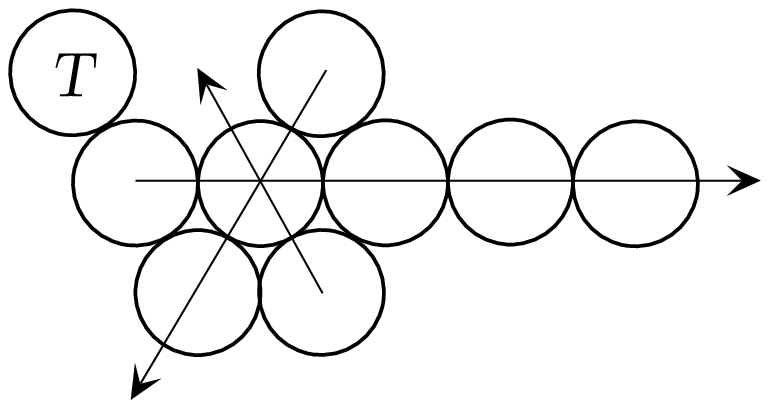}{.35}} 
\place{7}{15}{$\longrightarrow$} 
\kern 25pt
\lower10pt\hbox{\fig{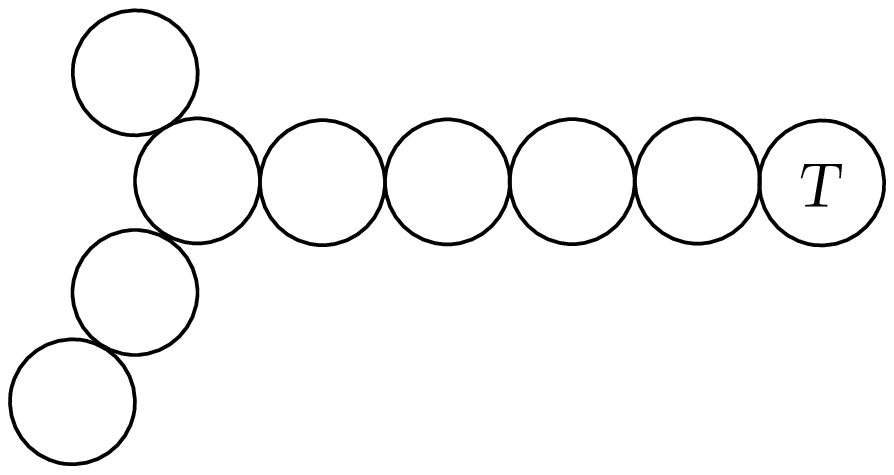}{.35}}}
\bigskip
\centerline{\small(a) Slips for II*}

\bigskip
\centerline{\fig{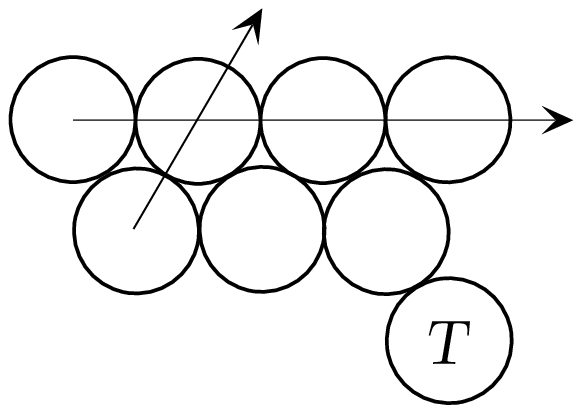}{.35}
\place{10}{15}{$\longrightarrow$}
\kern 40pt
\fig{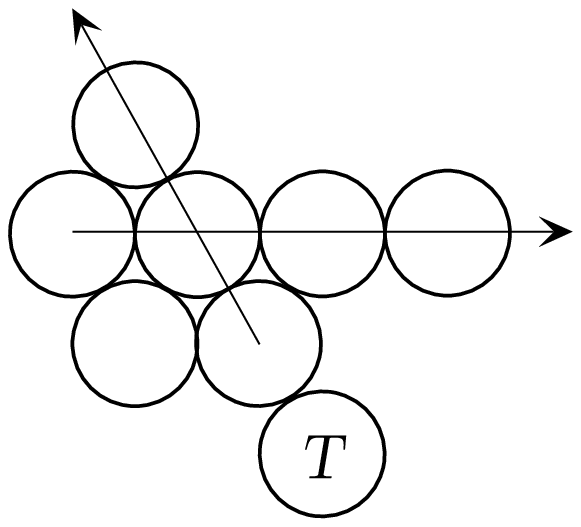}{.35} 
\place{10}{15}{$\longrightarrow$} 
\kern 40pt
\lower5pt\hbox{\fig{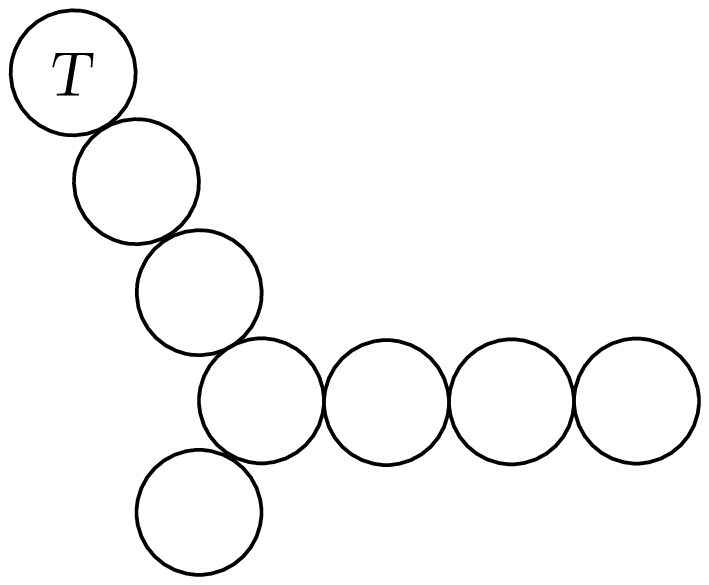}{.35}}}
\bigskip
\centerline{\small(b) Slips for III*}

\bigskip
\centerline{\fig{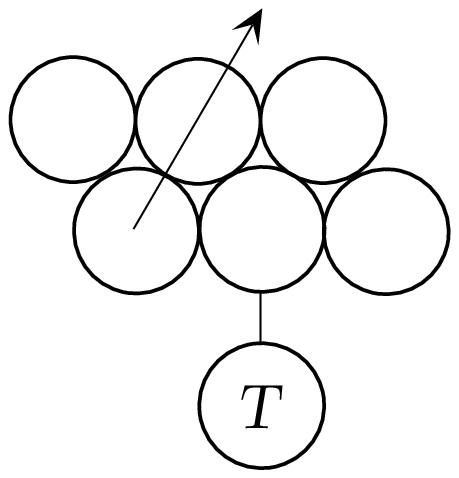}{.35}
\place{10}{20}{$\longrightarrow$}
\kern 40pt
\fig{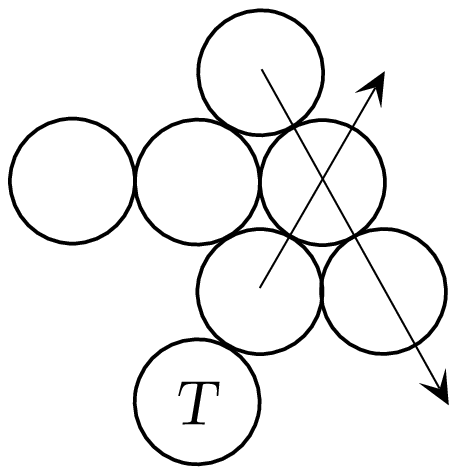}{.35} 
\place{10}{20}{$\longrightarrow$} 
\kern 40pt
\lower5pt\hbox{\fig{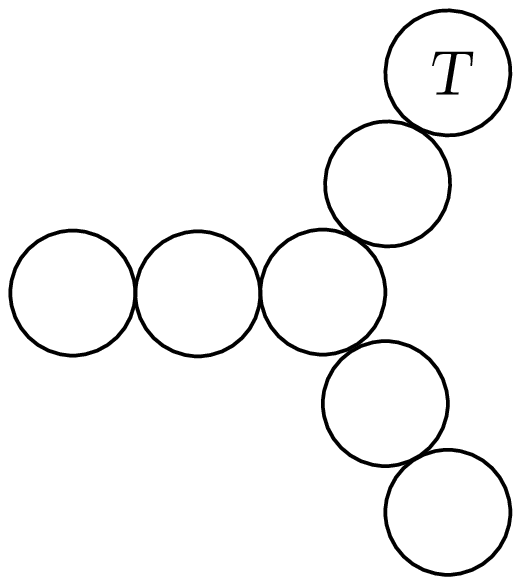}{.35}}} 
\bigskip
\centerline{\small(c) Slips for IV*}
\figno{2.7}
\end{figure}

\noindent
This completes the proof of Theorem 1*. 
\end{proof}


\medskip
\section{Gromov's compactness theorem}

In this section we state Gromov's compactness theorem for pseudo-holomorphic curves
\cite{Gro}, and then show how this theorem is richly illustrated by the convergence
of a sequence of regular fibers to a singular fiber in an elliptic surface.

\head{Pseudo-holomorphic curves and cusp curves}

A {\it pseudo-holomorphic curve} in an almost complex manifold $V$ is a smooth
map $f\co S\to V$ of a Riemann surface $S$ into $V$ whose differential at each point
is complex linear.  This means $df\circ j=J\circ df$, where $J$ is the almost
complex structure on $V$ (a bundle map on  $\tau_V$ with $J^2=-I$ on each fiber) and
$j$ is the (almost) complex structure on $S$.\foot{The complex structure on
$S$ is determined by $j$, since almost complex structures are integrable in
complex dimension $1$.}  We also write
$$
f\co (S,j)\to V
$$
to highlight the complex structure on $S$, which may vary; the almost
complex structure on $V$ is assumed fixed.   

Cusp curves are generalizations of pseudo-holomorphic curves in which one
allows the domain to be a {\it singular\/} Riemann surface $\bar S$,
obtained by crushing each component $C_i$ of a smoothly embedded 1--manifold
$C$ in $S$ to a point $p_i$.  Each singular point $p_i$ is to be viewed as a
{\sl transverse\/} double point of $\bar S$ with distinct complex
structures on the two intersecting sheets.  To make this precise, consider
(following Hummel \cite{Hum}) the {\sl smooth\/} surface $\hat S$ obtained
from $S-C$ by one-point compactifying each end separately, and let
$\alpha\co \hat{S} \rightarrow \overline{S}$ be the natural projection.  A
complex structure $\bar j$ on $\bar S$ is by definition a complex structure
$\hat j$ on $\hat S$.  The pair $(\bar S,\bar j)$ is then called a {\it
singular Riemann surface}, and a map
$$
\bar f\co (\bar S,\bar j) \to V
$$
is called a {\it cusp~curve\/} if $\hat f = \bar f \circ \alpha\co (\hat S,\hat
j) \to V$ is pseudo-holomorphic.  This setup is illustrated in Figure 3.1.
Note that because of dimension limitations a tangency has been drawn at the
singular points $p_i$, but these should be thought of as transverse
intersections; indeed $\bar f$ can map these to transverse double points as
for example occur in the core of a plumbing.   Examples of cusp curves are
given below.


\begin{figure}[ht!]
\centerline{\fig{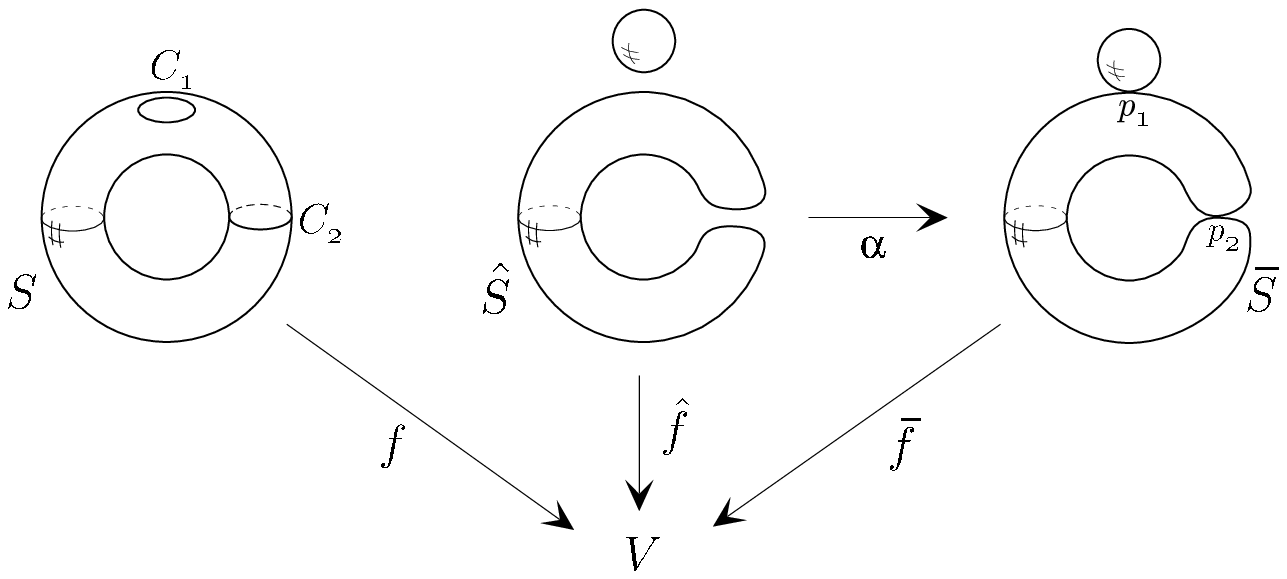}{.7}}
\figno{3.1}
\end{figure}

We also need the notion of a {\it deformation\/} of $S$ onto
$\bar S$, which by definition is any continuous map $d\co S\to \bar S$ which
sends each $C_i$ to $p_i$, and $S-\cup C_i$ diffeomorphically onto $\bar
S-\cup p_i$.

\head{Gromov's compactness theorem}

Let $f_k \co  (S,j_k) \rightarrow V$ (for $k = 1,2, \ldots$) be a sequence
of pseudo-holomorphic curves. Suppose that $V$ has a Hermitian metric and
that there is a uniform bound for the areas of the images $f_k(S)$.  (The
last condition is automatic if $V$ is symplectic and all $f_k(S)$ belong to
the same homology class \cite[page 82]{Hum}.)

Then Gromov's compactness theorem states that some subsequence of $f_k$
weakly converges to a cusp curve $\bar f\co  (\bar S,\bar j) \rightarrow V$,
where {\it weakly converges\/} means (for $k$ now indexing the subsequence):

\begin{enumerate}

\medskip
\item[(1)]  there exist deformations $d_k\co  S \rightarrow \bar S$ such that
the complex structures $(d_k^{-1})^*j_k$ on $\bar S$ minus the singular
points converge in the $C^{\infty}$ topology to $\bar j$ away from the
singular points;

\medskip
\item[(2)]  $f_k \circ d_k^{-1}$ converges in $C^{\infty}$ and
uniformly to $\bar f$ away from the singular points;

\medskip
\item[(3)]  the areas of the $f_k(S)$ converge to the ``area" of
$\bar f(\bar S)$ ( = area of $\hat f(\hat S)$).

\end{enumerate}

\head{Remarks}

\vskip-.1in
\begin{enumerate}

\medskip
\item[$\bullet$]  If the complex structures $j_k$ are all identical, then the
curves $C_i$ above must be null-homotopic and thus bound disks in $S$
(see \S3 in Chapter V of \cite{Hum}).  The standard example is given by
the sequence of pseudo-holomorphic curves $f_k \co  S^2 \rightarrow S^2 \times
S^2$  defined by $f_k (z) = (z, 1/(k^2z))$ which converges to the two curves
$S^2 \times 0 \cup 0 \times S^2$.  This process is called {\it bubbling off}
because of the appearance of the extra 2--sphere $0 \times S^2$.  A sequence
of circles that pinch to $(0,0)$ is $C_k = \{z:|z|=1/k\}$.  

\medskip
\item[$\bullet$]  When the complex structures $j_k$ are not identical,
then we may assume that they determine hyperbolic structures (taking out
three points if $S = S^2$ or one point if $S = T^2$) and then
$\bar S$ is the limit in the sense of degeneration of hyperbolic
structures. This will be the case with the fishtail singularities below.

\medskip
\item[$\bullet$]  That the deformations $d_k$ are necessary can be seen by
considering the case $S=V=T^2$ with $f_k\co T^2\to T^2$ equal to the $k^{\rm th}$
power of a Dehn twist around a meridian.  This sequence weakly converges to the
identity if we take $d_k = f_k$, but has no weakly convergent subsequence if
the $d_k$ are independent of $k$.

\end{enumerate}

\head{Examples}

Our collection of examples arise from the fact that a sequence of
regular fibers in an elliptic surface which converge to a singular fiber
form an illustration of Gromov's Theorem.  

To study this convergence, it is useful to have explicit models for the
projection of the elliptic surface near the singular fibers.  In particular,
neighborhoods of the singular fibers with finite monodromy --- namely those
of type I$^*_0$, II, II*, III, III*, IV and IV* --- can be obtained by taking
the quotient of $T^2 \times B^2$ by a finite group action followed by
resolving singular points and, perhaps, blowing down.  Since the
finite group action preserves a natural complex structure on $T^2$, it
follows that all regular fibers in these neighborhoods have the same
complex structure, so the circles that are pinched in the compactness
theorem are all null-homotopic.  Thus in these cases, only bubbling off
occurs.  In the other cases --- namely I$_n$ and I$^*_n$ for $b\geq 1$ --- the
monodromy is infinite, and the complex structure differs from fiber to
fiber.  Thus pinching of essential circles is allowed, and in fact necessary
in these cases.

We now discuss each of the different types of singular fibers.

\medskip
\ub{Type  I}  A neighborhood of the fishtail $I_1$ is obtained
from $T^2 \times B^2$ by adding a 2--handle to a vanishing cycle, for
example a
meridian of the torus.  This indicates that there is no bubbling off, and
that, after removing a point from the torus to make it hyperbolic, a
shortest geodesic representing the meridian shrinks (neck stretching) to a
point in the limit.  Near this point, using local coordinates $(z,w)$, the
projection to $S^2$ is given by the product $zw$ and the preimage of zero is
the two axes. Nearby the preimage is $zw = \epsilon$ which is an annulus. 
Note that the point that was removed to get a hyperbolic structure turns out
to be a removable singularity, that is, the pseudo-holomorphic map on the
punctured torus extends to the torus, and this is true in the limit.

For I$_n$ with $n >1$, there are $n$ parallel vanishing cycles, so we remove
$n$ points from $T^2$, interspersed so that the $n$ meridians are not
homotopic.  Then these meridians shrink as in the case of $I_1$, and we get a
necklace of $n$ spheres in the limit.  Note that these spheres have
multiplicity one, so their sum is homologous to a regular fiber and hence
has square zero. This implies that each 2--sphere has self-intersection $-2$.

\medskip
\ub{Type  I*}  First consider the case I$^*_0$.  The monodromy is of
finite order, namely two, so we begin with a simple $\Z_2$--action.  View the
torus as the unit square with opposite sides identified, or equivalently the
quotient of $\C$ by the lattice $\langle1,i\rangle$ (with the induced
complex structure).  Let $\sigma_2$ be the involution on $T^2$ which rotates
the square (or $\C$) by $\pi$ (Figure\ 3.2a).  Clearly the quotient
$X=T^2/\sigma_2$ is a 2--sphere, and the projection $T^2\to X$ is a 2--fold
cover branched over four points corresponding to the center $a$ and vertex $b$ of
the square, and the two pairs $c$ and $c'$ of midpoints on opposite edges.


\begin{figure}[ht!]
\centerline{\fig{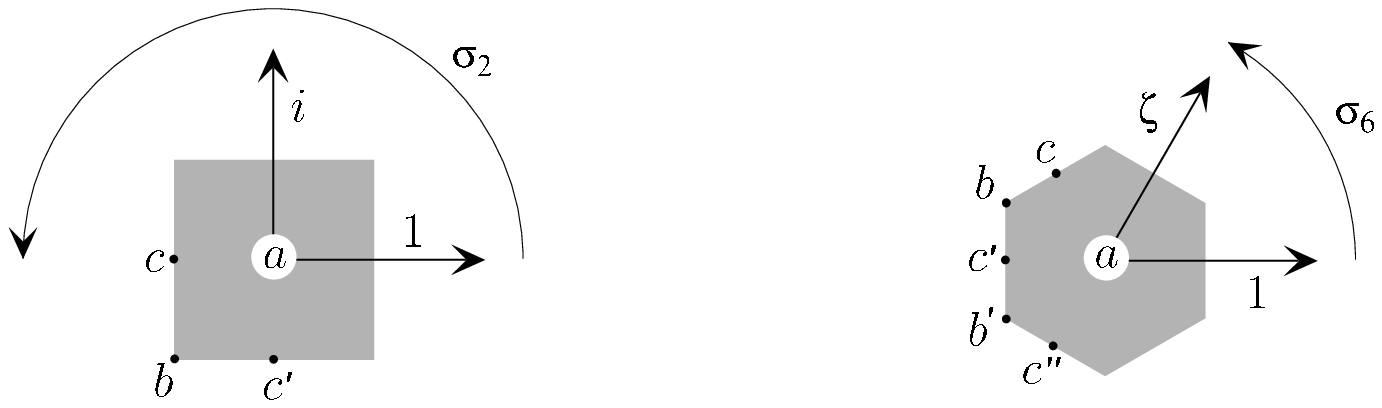}{.7}}
\smallskip
\centerline{\small(a) \kern 155pt (b) \kern5pt}
\figno{3.2: Automorphisms of $T^2$}
\end{figure}

Set $E = T^2 \times B^2 / \sigma_2 \times
\tau_2$ and $D = B^2/ \tau_2 \cong B^2$, where $\tau_r$ is the rotation of
$B^2$ by $2\pi /r$.  Then the natural projection $E \to D$ has a singular
fiber over $0\in D$, namely $X$, and is a $T^2$--bundle over $D - 0$ with
monodromy $-I$.  The space $E$ is a 4--manifold except at the four branch
points on $X$ which are locally cones on $\R P^3$.  These singular points
can be resolved by removing the open cones on $\R P^3$ and gluing in
cotangent disk bundles of $S^2$ with cores $A,\ B,\ C$ and $D$.  This gives
a neighborhood $N($I$^*_0)$ of the singular fiber I$^*_0$, as shown in
Figure 3.3.  Note that $2X + A + B + C + D$ is homologous to a regular fiber
and hence has square zero, which implies that $X\cdot X = -2$.  The second
picture is the standard (dual) plumbing diagram, where the vertex weights
are the multiplicities.


\begin{figure}[ht!] 
\centerline{\fig{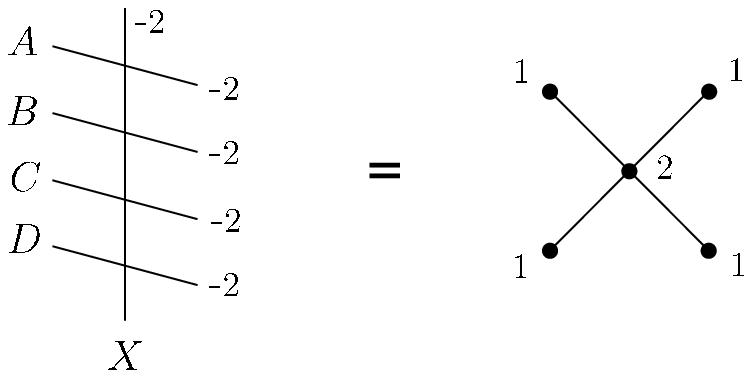}{.7}}
\figno{3.3: $N($I$_0^*)$}
\end{figure}

Now the picture for Gromov's compactness theorem is clear (Figure\ 3.4).  The
torus bubbles off four 2--spheres at the branch points (labelled $a$--$d$)
and then double covers $X$ while the four bubbles hit $A$--$D$ with
multiplicity one.  (Note that the complex structures on all the torus fibers
are the same since the group action is holomorphic, and so we expect to see
only bubbling off.)  The picture is drawn so that the branched covering
transformation $\sigma_2$ corresponds to a $\pi$--rotation about the vertical
axis through the branched points.  As above, tangencies in the picture
of I$_0^*$ correspond to transverse double points in $N($I$_0^*)$.


\begin{figure}[ht!] 
\centerline{\fig{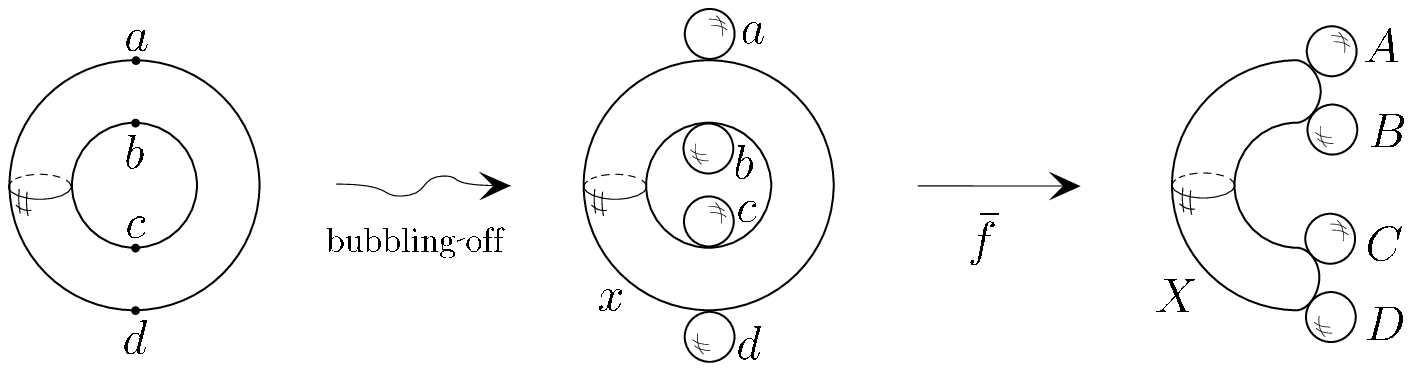}{.7}}
\smallskip
\centerline{\small$S=T^2$ \kern 90pt $\bar S$ \kern 110pt I$_0^*$ \kern10pt}
\figno{3.4: Degeneration to I$_0^*$}
\end{figure}

For I$^*_n$ with $n>0$, as in the cases I$_n$ above, the monodromy is of
infinite order.  Nonetheless, this case is related to I$^*_0$ in that we
perform the same construction and in addition shrink $n$ pairs of meridional
circles, each pair having the same image under the double covering map.  We
get $\bar{S}$ as drawn in Figure 3.5, and it is mapped onto the singular
fiber in the obvious way: each sphere labelled with a lower case letter maps
one-to-one onto the sphere with the corresponding upper case letter and
subscript, with the exception of $x_1$ and $x_n$ which map by 2--fold covers
to $X_1$ and $X_n$.


\begin{figure}[ht!] 
\centerline{\fig{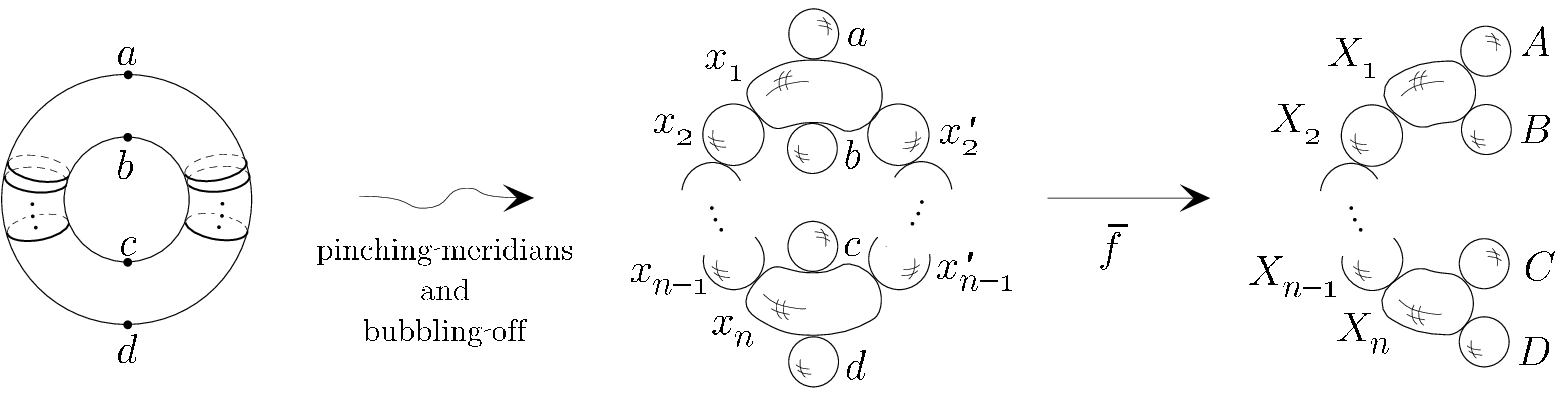}{.8}}
\smallskip
\centerline{\small$S=T^2$ \kern 125pt $\bar S$ \kern 135pt I$_n^*$}
\figno{3.5: Degeneration to I$_n^*$}
\end{figure}

\medskip
\ub{Type  II}  
View the torus as the hexagon with opposite sides identified, or
equivalently the quotient $\C/\langle 1,\zeta \rangle$ where $\zeta =
\exp(2\pi i/6)$.  Let $\sigma_6$ be the automorphism of $T^2$ of order 6
which rotates the hexagon (or $\C$) by $2\pi/6$.  The generic orbit of this
action has six points. However, the center $a$ of the hexagon is a fixed
point with stabilizer $\Z_6$, the vertices form an orbit with two points
$b,b'$ and stabilizer $\Z_3$, and the midpoints of the sides form an orbit
with three points $c,c',c''$ and stabilizer $\Z_2$ (see Figure\ 3.2b).  The
quotient $X=T^2/\sigma_6$ is again a 2--sphere, and the projection $T^2\to X$
is a 6--fold (irregular) cover branched over three points with branching
indices 6, 3, and 2.

Set $E = T^2 \times B^2 / \sigma_6 \times \tau_6$ and $D = B^2/
\tau_6 \cong B^2$.  The projection $E \to D$ has the singular 2--sphere
$X$ over $0 \in D$, and is a bundle over $D - 0$ with monodromy ${\ 1\
\,1\,\choose-1\ 0\,}$ ($=UV$ for $U$ and $V$ as in \S2).  As before, $E$
is a manifold except at the three branch points on $X$ which are locally
cones on $L(6,1)$, $L(3,1)$, and $L(2,1)$.  These singularities can be
resolved by cutting out the cones and gluing in the disk bundles over $S^2$
of Euler class $-6$, $-3$, and $-2$ respectively.  The torus fiber is
homologous to $6X + A + 2B + 3C$, and setting its square equal to zero and
solving gives $X\cdot X = -1$.  This gives the first picture in Figure 3.6,
which is followed by a sequence of blowdowns to produce the neighborhood
$N($II) of the cusp.


\begin{figure}[ht!] 
\centerline{\fig{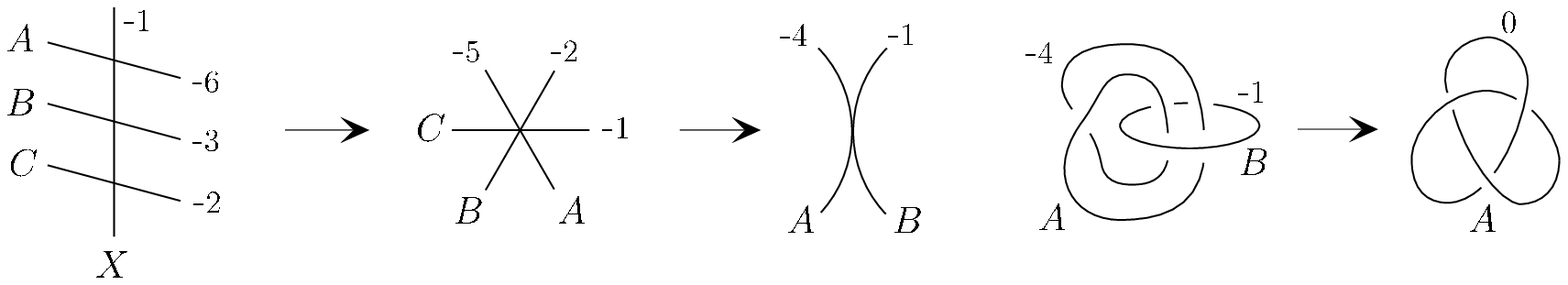}{.7} \place{-140}{35}{$=$}}
\figno{3.6: Blowing down to $N($II)}
\end{figure}

In the compactness theorem, the torus fiber bubbles off a 2--sphere at
each of the six branch points $a,b,b',c,c',c''$.  The torus then 6--fold
covers $X$ while bubble $a$ hits $A$, the two bubbles $b,b'$ hit $B$, and
the three bubbles $c,c',c''$ hit $C$.  Finally this map is composed with the
sequence of three blowdowns to give a degeneration to the cusp.  Since the
blowdowns burst all but the first bubble, this amounts to one bubble
which maps onto the cusp while the torus is mapped (by the constant
holomorphic map) to the singular point of the cusp (Figure\ 3.7).


\begin{figure}[ht!]
\centerline{\fig{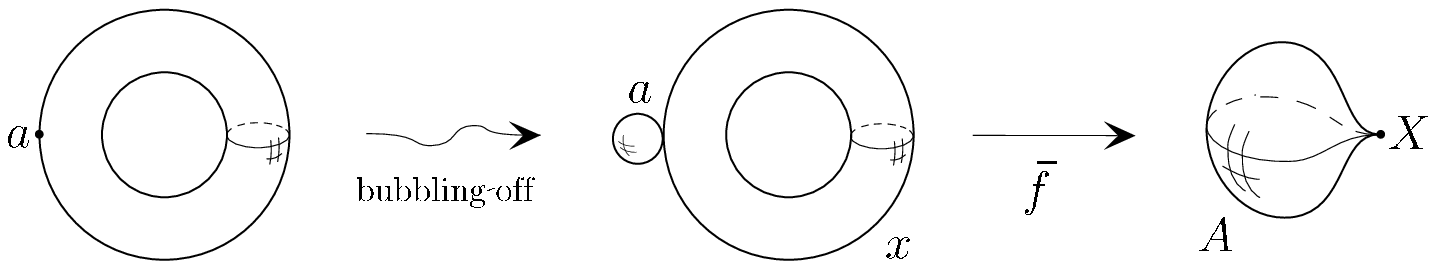}{.75}}
\centerline{\small $S=T^2$ \kern 100pt $\bar S$ \kern 100pt II \kern10pt}
\figno{3.7: Degeneration to II}
\end{figure}

\medskip
\ub{Type  II*}    We use the same construction for $E$ as in the
previous case for II, but the orientation is changed.  This is because
the rational elliptic surface $E(1)$ equals II and II* glued along their
common boundary.  So $E$ is now desingularized by removing the cones on
$L(6,5)$, $L(3,2)$, and $L(2,1)$ and replacing them by linear plumbings (of
disk bundles of Euler class $-2$) of length 5, 2 and 1 respectively (Figure\
3.8).  This is $N$(II*); note that $X \cdot X = -2$ by the usual
calculation.


\begin{figure}[ht!] 
\centerline{\fig{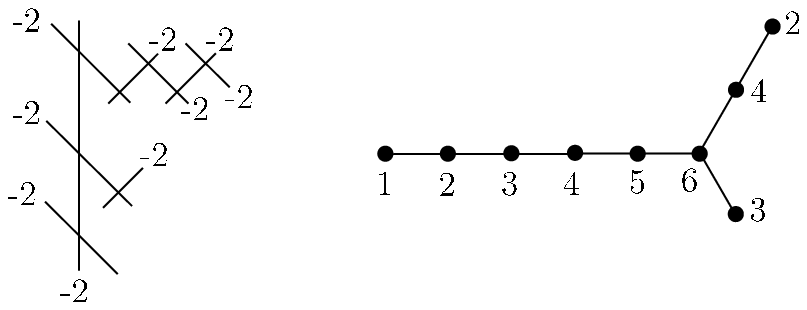}{.8} \place{-137}{37}{$=$}}
\figno{3.8: $N($II*)}
\end{figure}

For the compactness theorem, the torus bubbles off a linear graph of five
bubbles at the fixed point $a$, a line of two bubbles at each of the two
points $b,b'$ with stabilizer $\Z_3$, and one bubble at each of the three
points $c,c',c''$ with stabilizer $\Z_2$ (see Figure 3.9).  Then the torus
6--fold covers $X$; each of the five bubbles in the linear graph 5--fold cover
(with two branch points), 4--fold cover, 3--fold cover, 2--fold cover (still
with two branch points where they intersect their neighbors), and 1--fold
cover, the long arm of II*; the pair of two bubbles each will 2--fold cover
and 1--fold cover, providing multiplicities 4 and 2 since there are two
pairs; the three single bubbles all map onto the short arm of II* giving
multiplicity 3.  Note that the labels for the bubbles have been chosen so
that bubble $a_i$ maps to $A_i$ by an $i$--fold branched covering, and
similarly for the $b$ and $c$--bubbles.


\begin{figure}[ht!] 
\centerline{\fig{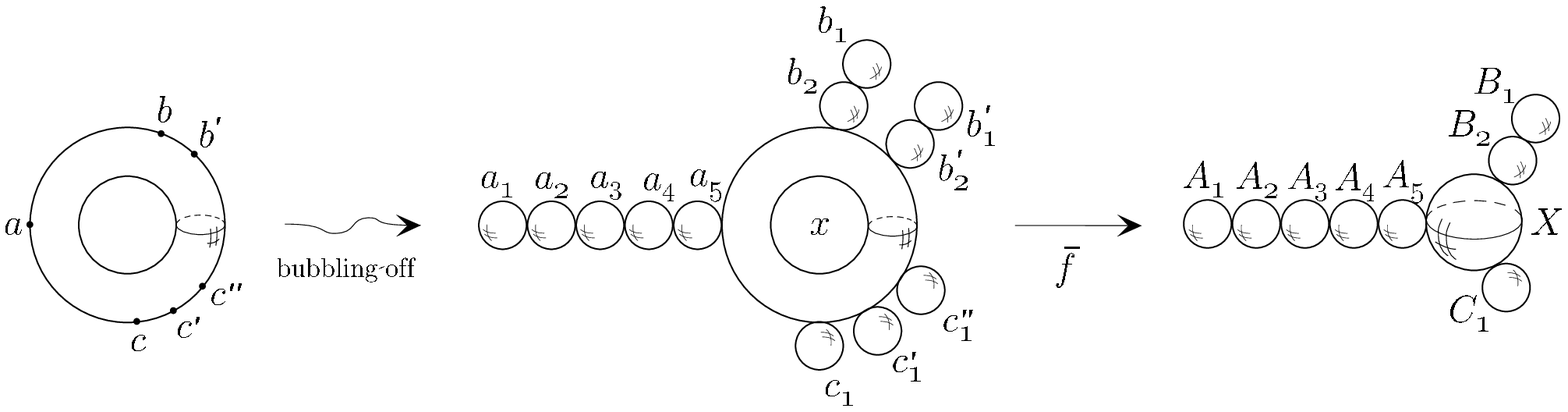}{.65}}
\centerline{\small$S=T^2$ \kern 125pt $\bar S$ \kern 135pt II*}
\figno{3.9: Degeneration to II*}
\end{figure}

We can now abbreviate the description for III and IV and their
duals for the arguments are similar to II and II* with no new
techniques.

\medskip
\ub{Type  III}    Again consider the torus as the square with
opposite sides identified and let $\sigma_4$ be rotation by $\pi /2$. This
has fixed points at the center of the square and at the vertex, and an
orbit of two points equal to the midpoints of the sides with stabilizer
$\Z_2$.  We resolve the quotient $T^2 \times B^2/ \sigma_4 \times \tau_4$ by
cutting out cones and gluing in disk bundles, and a sequence of blowdowns gives
the neighborhood $N($III) shown in Figure 3.10.


\begin{figure}[ht!] 
\centerline{\fig{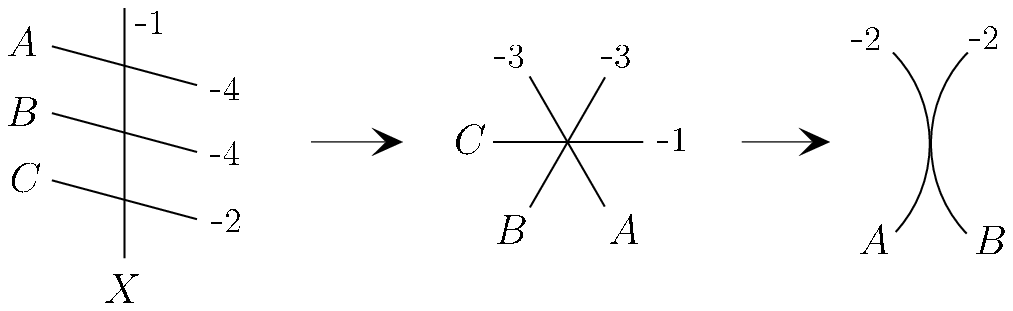}{.8}}
\figno{3.10: Blowing down to $N($III)}
\end{figure}

Now, as in II, the torus bubbles off 2--spheres which map and then in
some cases blow down, to give a composition in which two bubbles
survive and map onto the two curves in III and the torus maps to the
point of tangency (Figure\ 3.11).


\begin{figure}[ht!] 
\centerline{\fig{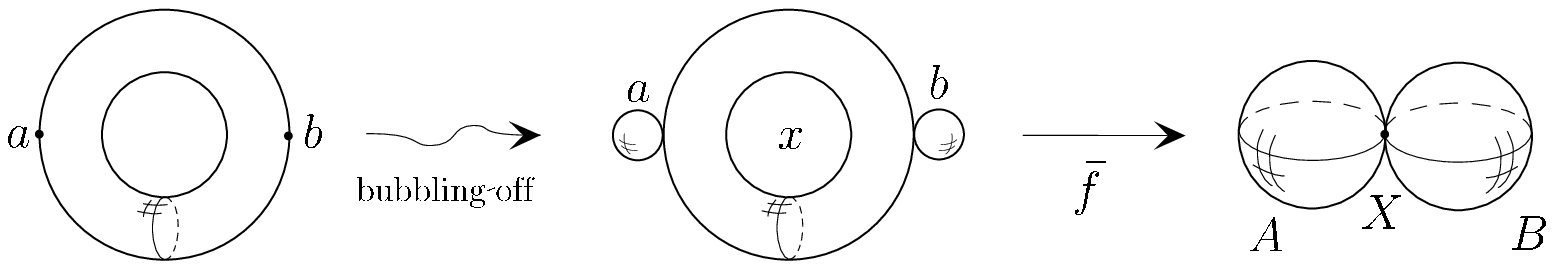}{.75}}
\centerline{\small$S=T^2$ \kern 115pt $\bar S$ \kern 115pt III \kern10pt}
\figno{3.11: Degeneration to III}
\end{figure}

\medskip
\ub{Type  III*}    As with II*, we reverse orientation, cut out
the cones and replace them with linear plumbings to get $N$(III*) (Figure\
3.12).


\begin{figure}[ht!]
\centerline{\fig{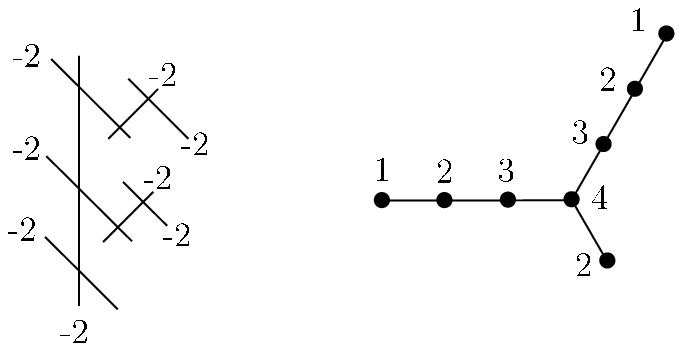}{.8} \place{-110}{37}{$=$}}
\figno{3.12: $N($III*)}
\end{figure}

We then get a degeneration of the torus fiber very similar to II*, using the
torus as a $4$--fold branched cover of $S^2$ with three branch points of
indices $4$, $4$, and $2$ (Figure\ 3.13).


\begin{figure}[ht!] 
\centerline{\fig{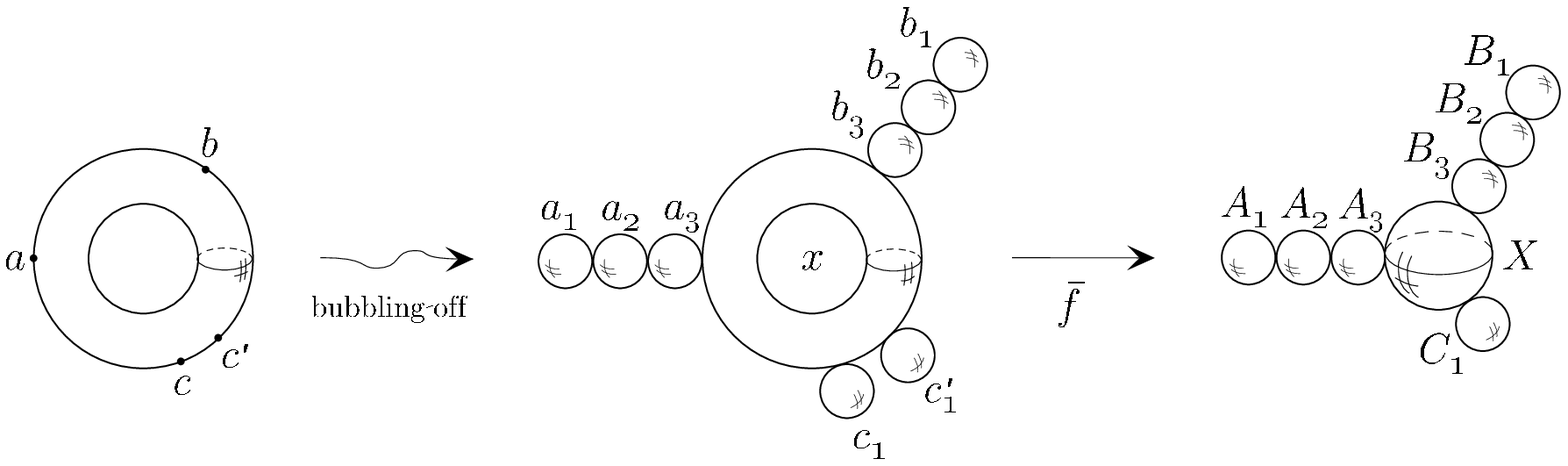}{.7}}
\centerline{\small$S=T^2$ \kern 125pt $\bar S$ \kern 135pt III*}
\figno{3.13: Degeneration to III*}
\end{figure}

\ub{Type  IV}    Here we define a $\Z_3$ action on the torus by
simply squaring the action given in II.  The center and any two adjacent
vertices of the hexagon represent the three fixed points $a$, $b$ and $c$
($\,= b'$ in Figure\ 3.2b).  Proceeding as before, we get the singular fiber
drawn in Figure 3.14.  Blowing down once gives
$N$(IV).


\begin{figure}[ht!] 
\centerline{\fig{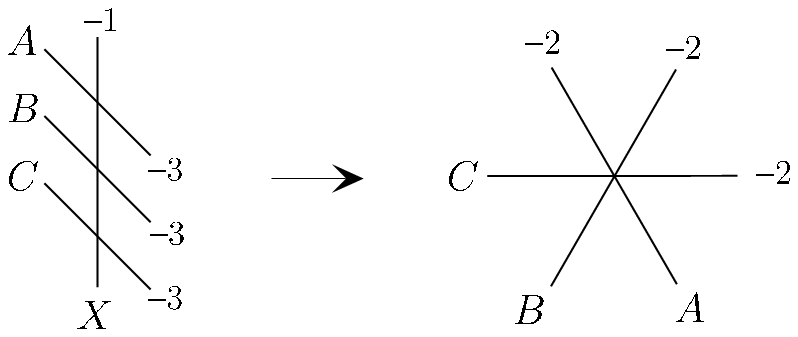}{.8}}
\figno{3.14: Blowing down to $N($IV)}
\end{figure}

The torus now bubbles off three 2--spheres which hit the three curves in IV
while the torus maps to the triple point (Figure\ 3.15).


\begin{figure}[ht!] 
\centerline{\fig{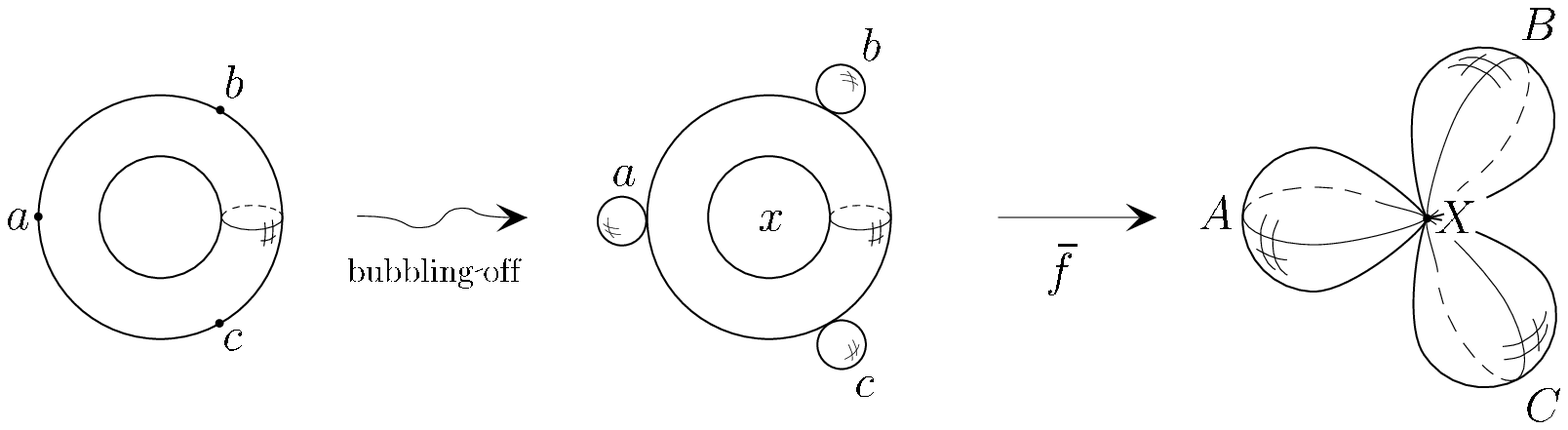}{.75}}
\centerline{\small$S=T^2$ \kern 120pt $\bar S$ \kern 140pt IV}
\figno{3.15: Degeneration to IV}
\end{figure}

\medskip
\ub{Type  IV*}    Arguments similar to those above give $N$(IV*)
(Figure\ 3.16) as well as the degeneration of the torus fiber (Figure\ 3.17).


\begin{figure}[ht!] 
\centerline{\fig{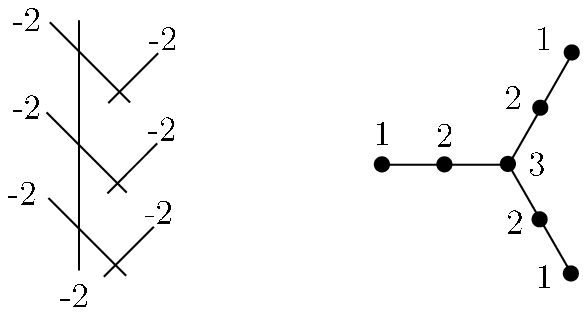}{.8} \place{-90}{37}{$=$}}
\figno{3.16: $N($IV*)}
\end{figure}


\begin{figure}[ht!] 
\centerline{\fig{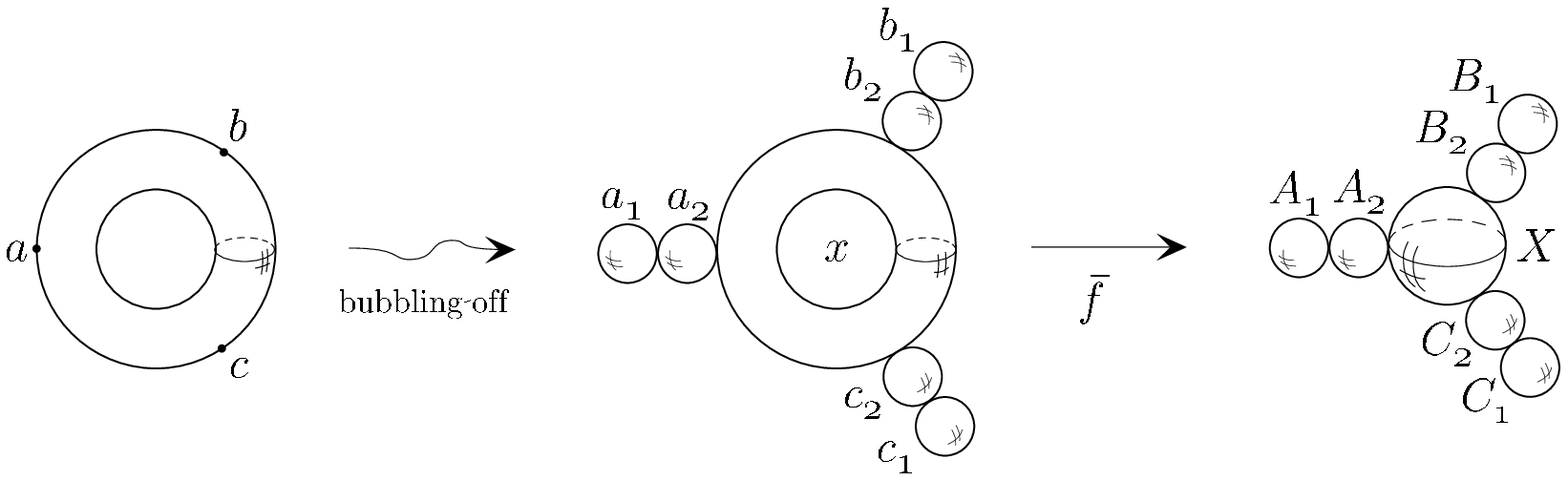}{.75}}
\centerline{\small$S=T^2$ \kern 130pt $\bar S$ \kern 120pt IV* \kern10pt}
\figno{3.17: Degeneration to IV*}
\end{figure}


\Addresses\recd

\end{document}